\documentclass[11pt]{amsart}
\usepackage{amssymb,latexsym}
\newtheorem{theorem}{Theorem}[section]

\newtheorem{lemma}[theorem]{Lemma}
\newtheorem{remark}[theorem]{Remark}

\newtheorem{definition}[theorem]{Definition}
\newtheorem{cor}[theorem]{Corollary}

\begin{document}

\title{Birational cobordism invariance of uniruled symplectic manifolds}

\author{Jianxun Hu$^1$ \& Tian-Jun Li$^2$ \& Yongbin Ruan$^3$}
\address{Department of Mathematics\\ Zhongshan
University\\Guangzhou, 510275 \\ P. R. China}
\email{stsjxhu@mail.sysu.edu.cn}
\thanks{${}^1$Partially supported by the NSFC Grant 10631050, NKBRPC(2006CB805905) and NCET-04-0795}
\address{School  of Mathematics\\  University of Minnesota\\ Minneapolis, MN 55455}
\email{tjli@math.umn.edu}
\thanks{${}^2$supported by NSF Grant}
\address{Department of Mathematics\\ University of Michigan\\ Ann Arbor, MI
48109-1109\\  and \\ Yangtze Center of Mathematics \\ Sichuan
University \\ Chengdu, 610064, P.R. China
 }\email{ruan@umich.edu}
 \thanks{${}^3$supported by NSF Grant}

\maketitle

\tableofcontents

\section{Introduction}

Birational geometry has always been an important topic in algebraic
geometry. In the 80's, an industry called Mori's birational geometry
program was created for the birational classification of algebraic
manifolds of dimension three.
 In the early 90's, the last
author observed that some aspects of this extremely rich program of
Mori can be extended to symplectic geometry via the newly created
Gromov-Witten theory \cite{R1}. Further, he speculated that in fact
there should be a symplectic birational geometric program. Such a
program is important in two ways. The flexibility of symplectic
geometry should give a better understanding of birational algebraic
geometry in the same way that the Gromov-Witten theory gave a much
better understanding of the role of rational curves in Mori theory.
Secondly, such a symplectic birational geometry should be the first
step towards a classification of symplectic manifolds \cite{R2}.
During the last ten years, there was virtually no progress in this
direction. There are many reasons for lack of this progress, one  of
which  was the difficulty to carry out computations in Gromov-Witten
theory. Fortunately, a great deal of progress has been made to
remedy this aspect of the problems. Many techniques have been
developed to calculate the Gromov-Witten invariants. It seems to be
a good time now  to restart a real push for the symplectic
birational geometry. This is the first of a series of papers by the
last two authors to treat this new subject of symplectic birational
geometry. Our treatment is by no means complete. On the contrary,
there are many more problems being discovered than  answered.

For a long time, it was not really clear  what is an appropriate
notion of birational equivalence in symplectic geometry. Simple
birational operations such as blow-up/blow-down were known in
symplectic geometry for a long time \cite{GS, MS}. But there is no
straightforward generalization of  the notion of a general
birational map in the flexible symplectic category. The situation
changed a great deal when the weak factorization theorem was
established recently (see the lecture notes \cite{M} and the
reference therein) that any birational map between projective
manifolds can be decomposed as a sequence of blow-ups and
blow-downs. This fundamental result resonates perfectly with the
picture of the wall crossing of symplectic reductions analyzed by
Guillemin-Sternberg in the 80's. Therefore, we propose to use their
notion of cobordism in \cite{GS} as the symplectic analogue of the
birational equivalence (see Definition \ref{bc}). To avoid confusion
with other notions of cobordism in the symplectic category, we would
call it {\it symplectic birational cobordism}.

A fundamental concept of birational geometry is uniruledness.
Algebro-geometrically, it means that the manifold is covered by
rational curves. Notice that, by \cite{Ltj}, it is not meaningful to
define this notion by simply mimicking the definition in algebraic
geometry and requiring that there is a symplectic sphere in a fixed
class through each point. Otherwise, every simply connected manifold
would be uniruled. On the other hand, by a theorem of Kollar-Ruan
\cite {K}, \cite{R1}, a uniruled projective manifold has a nonzero
genus zero GW-invariant with a point insertion. Therefore, we call a
symplectic manifold $(M,\omega)$ {\em (symplectically) uniruled} if
there is a nonzero genus zero GW invariant involving a point
constraint. Then, it is a fundamental problem in symplectic
birational geometry to prove that symplectic uniruledness is a
birational invariant. It is obviously deformation invariant. The
main purpose of this paper is to prove

\begin{theorem}
Symplectic uniruledness is invariant under symplectic blow-up and
blow-down.
\end{theorem}

 This theorem follows from  a general Relative/Absolute
 correspondence for a symplectic manifold together with
 a symplectic submanifold. Such a correspondence was first
 established in \cite{MP} when the submanifold is of codimension 2,
 i.e. a symplectic divisor.

Clearly our theorem can also be viewed  as a kind of blow-up/down
formula of GW invariants. It is rare to be able to obtain a
general blow-up formula. For the last ten years, only a few
limited cases were known \cite{H1, H2, H3, HZ, Ga}. Although our
technique is more powerful, for instance, it still does not work
for rational connectivity. Here a symplectic manifold is said to
be {\em rationally connected} if there is a nonzero genus zero GW
invariant involving two point insertions. We do speculate that
{\em symplectic rational connectivity is invariant under blow-up
and blow-down.}

A natural question  is whether Hamiltonian $S^1$ manifolds are
uniruled. Based on our blowup technique as well as Seidel
representation of $\pi_1$ of the Hamiltonian group in the small
quantum homology, McDuff \cite{Mc3} gave an affirmative answer to
this question.


For the reader's convenience, we review in sections two and three
some background materials scattered in the literature. In section
two, we will review the definition of birational cobordism in the
symplectic geometry and describe blow-up, blow-down and $\mathbb
Z-$linear deformation as birational cobordisms. In section three,
we will briefly review the relative GW-invariant and the
degeneration formula which is an essential tool for us. In section
four, we will give the definition of symplectic uniruledness and
some elementary properties of uniruled symplectic manifolds. The
core of the paper consists of  sections five and six. In section
five, we will describe and partially verify the Relative versus
Absolute Correspondence for GW-invariants, which generalizes early
works of Maulik-Okounkov-Pandharipande \cite{OP}, \cite{MP}. Then,
in section six, we apply our Relative versus Absolute
correspondence to prove the main theorem. In section seven, we
will complete the proof of the Relative versus Absolute
correspondence by computing certain relative GW invariants of
$({\mathbb P}^n, {\mathbb P}^{n-1})$.

 The first author would like to thank Z. Jiang and B. Zhang for
 their help. The third  author would like to thank Y.P. Lee for valuable discussions.
 We thank V. Guillemin for his interest and  J. Dorfmeister for
carefully reading this manuscript and improving its presentation. We
also thank the referees for pointing out vague statements and ways
to clarify them.
 Finally we are very grateful to D. McDuff for
many constructive comments and suggestions.

\section{Birational cobordism}

The basic reference for this section is \cite{GS}. We start with the
definition which is essentially contained in \cite{GS}.

\begin{definition} \label{bc} Two symplectic manifolds $(X, \omega)$ and
$(X', \omega')$ are birational cobordant if there are a finite
number of symplectic manifolds $(X_i, \omega_i)$,$ 0\leq i\leq k$,
with $(X_0, \omega_0)=(X, \omega)$ and $(X_k, \omega_k)=(X',
\omega')$, and for each $i$, $(X_i, \omega_i)$ and $(X_{i+1},
\omega_{i+1})$ are symplectic reductions of a semi-free Hamiltonian
$S^1$ symplectic manifold $W_i$ (of 2 more dimension).

\end{definition}

Here an $S^1$ action is called semi-free if it is free away from the
fixed point set.

There is  a related notion in dimension 4 in \cite{OO}. However we
remark that the cobordism relation studied in this paper is  quite
different from some other notions of symplectic cobordisms, see
\cite{EGH}, \cite{EH}, \cite{Gi}, \cite{GGK}.

According to \cite{GS}, we have the following basic factorization
result.

\begin{theorem} \label{f}
 A birational cobordism can be decomposed
as a sequence of elementary ones, which are modeled on blow-up,
blow-down and $\mathbb Z-$linear deformation of symplectic
structure.
\end{theorem}

Comparing with the weak factorization theorem, we then have

\begin{theorem}
 Two birational projective manifolds
with any polarizations are birational as symplectic manifolds.
\end{theorem}

We will now review each of these elementary birational cobordisms in
Theorem \ref{f}.

\subsection{Coupling form and linear deformations} \label{universal construction}
\subsubsection{Universal construction} \label{coupling form}
Let us first review the Sternberg-Weinstein universal construction.
Let $\pi:P\to X$ be a principal bundle with structure group $G$ over
a symplectic manifold $X$ with symplectic form $\omega$. If ${\rm
g}$ denote the Lie algebra of $G$, then  a connection on $P$ gives
rise to a $\rm g$ valued $1-$form on $P$ corresponding to the
projection onto the vertical. Let Ver be the vertical bundle of the
fibration. A $G-$invariant complementary subbundle $F$ is nothing
but a connection of $P$. It also embeds $P\times {\rm g}^*$ into
$T^*P$.

The desired $1-$form at $(p, \tau)$ is given by $\tau\cdot A$, where
we use $\cdot$ to denote the pairing between $\rm g$ and $\rm g^*$.
Denote this $1-$form by $\tau\cdot A$ as well. Notice that it is the
restriction of the canonical $1-$form on $T^*P$. Therefore
$d(\tau\cdot A)$ is non-degenerate on the fibers of $P\times {\rm
g}^*$.

  Then $\omega_A=\pi^*\omega+d(\tau\cdot A)$ is sometimes
called the coupling form of $A$. The $G-$action on $P\times \rm g^*$
given by
$$g(p, \zeta)=(g^{-1}p, \hbox{Ad}(g)^*\zeta),$$
preserves $\tau\cdot A$ and hence $\omega_A$.

Notice that at any point in $P\times \{0\}$, $\tau=0$ and
$\omega_A$ is equal to $\pi^*\omega+d\tau \cdot A$, hence it is
symplectic there.

\begin{lemma} \label{W} $\omega_A$ is
a symplectic form on $P\times {\mathcal W_A}$ for some $G-$invariant
neighborhood ${\mathcal W_A}$ of $0\in \rm g^*$. The projection onto
${\mathcal W_A}$ is a moment map  on $P\times {\mathcal W_A}$.
\end{lemma}

This lemma is well-known. Notice that the vertical bundle Ver is the
bundle of the null vectors of $\pi^*\omega$. So as explained in
\cite{GS}, the construction is just a special case of the
coisotropic embedding theorem. It follows from the uniqueness part
of the coisotropic embedding theorem that the symplectic structure
$\omega_A$ on $P\times {\mathcal W}_A$ near $P\times 0$ is
independent of $A$ up to symplectomorphisms.

More generally, if $(F, \omega_F)$ is a symplectic manifold with a
Hamiltonian $G$ action, we can form the associated   bundle
$P_F=P\times_G F$. Let $\mu_F:F\to \rm g^*$ be a moment map.
Furthermore, assume that
\begin{equation} \label{W2}
\mu_F(F)\subset \mathcal W_A.
\end{equation}
Then there is a symplectic structure $\omega_{F, A}$ on $P_F$ which
restricts to $\omega_F$ on each fiber.

To construct  $\omega_{F, A}$ consider the $2-$form
$\omega_A+\omega_F$ on $P\times \rm g^*\times F$. It is invariant
under the diagonal $G-$action and is symplectic on $P\times
{\mathcal W_A}\times F$.

The $G-$action is in fact Hamiltonian with
$$\Gamma_{\mathcal W_A}=\pi_{\rm g^*}+\mu_F:P\times {\mathcal W_A}\times F\to g^*$$
as a moment map.

Furthermore, by (\ref{W2}), for any $f\in F$, we have $\mu_F(f)\in
{\mathcal W_A}$, thus
$$\Gamma_{\mathcal W_A}^{-1}(0)=\{(p, -\mu_F(f), f)\}.$$
In particular, $\Gamma_{\mathcal W_A}^{-1}(0)$ is $G-$equivariantly
diffeomorphic to $P\times F$, and the symplectic reduction at $0$
yields the desired symplectic form $\omega_{F, A}$ on $P_F$.

 In fact when $(F, \omega_F)=(TG^*, \omega_{can})$,
 then  $P\times \rm g^*=P\times_G TG^*$.

\subsubsection{${\Bbb Z}-$linear deformation} \label{linear deformation}

\begin{definition}
A ${\mathbb Z}-$linear deformation is a path of symplectic form
$\omega+t\kappa$, $t\in I$, where $\kappa$ is a closed $2-$form
representing an integral class and $I$ is an interval. Two
symplectic forms are ${\Bbb Z}-$linear deformation equivalent if
they are joined by a finite number of ${\Bbb Z}-$linear
deformations.
\end{definition}

 Let $P$ be the principal
 $S^1-$bundle whose Chern class is $[\kappa]$. Let $A$ be a connection
 $1-$form such that $dA=\pi^*\kappa$.
Write $\omega_A=\pi^*\omega+d(tA)$, where $t$ is the linear
coordinate on  $\rm g^*={\mathbb R}^*$. It follows from Lemma
$\omega_A$ is symplectic at $(x, \theta, t)\in P\times {\mathbb
R}^*$ if and only if $\omega+t\kappa$ is symplectic at $x\in X$.
Thus the Duistermaat--Heckman Theorem can be interpreted as saying
that a
 ${\Bbb Z}-$linear deformation is a birational cobordism.


It seems to be more natural to consider a general deformation of
symplectic structures. However the following lemma shows that the
two notions are essentially the same.

\begin{lemma} Let $\omega_t, t\in [0,1],$ be a path of symplectic forms with $[\omega_0]-[\omega_1]$
being rational. Then $\omega_0$ and $\omega_1$ are ${\Bbb Z}-$linear
deformation equivalent.
\end{lemma}

\begin{proof} Observe that for any  convex open neighborhood of the space of symplectic forms
the convex linear combination of any two symplectic forms with
relative rational period is a ${\Bbb Z}-$linear deformation. Observe
also that any $\omega$ has a convex neighborhood.

 For each $x\in [0,1]$ let $I_x$ be a
neighborhood such that the family  $\omega_t, t\in I_x,$ is in a
convex neighborhood $V_{\omega_x}$. Let $I_{x_1},  \cdots , I_{x_n}$
with $0=x_0<  \cdots <x_n=1$ be a finite subcover of $[0,1]$. Since
$V_{\omega_{x_0}}$ and $V_{\omega_{x_1}}$ are both convex and
rational symplectic forms are dense in the nonempty set
$V_{\omega_{x_0}}\cap V_{\omega_{x_{1}}}$, $\omega_0$ is ${\Bbb
Z}-$linear deformation equivalent to any symplectic form $\omega\in
V_{\omega_{x_1}}$ with relative rational period. By repeating this
argument we find that $\omega_0$ is ${\Bbb Z}-$linear deformation
equivalent to $\omega_1$.

\end{proof}

\begin{remark} \label{deformation2} For any symplectic form $\omega$ and $J$ tamed by $\omega$,
there is a nearby $\omega'$ with rational period and still tamed by
$J$. Further, $\omega'$ also tames any $J'$ near $J$. Thus given any
path  of symplectic forms  $\omega_t, t\in [0,1]$ and $J_t$ be a
path of almost complex structures such that $J_t$ is tamed by
$\omega_t$ for each $t$, we can find finitely many symplectic forms
with rational period such that each $J_t$ is tamed by one of them.
\end{remark}

\subsection{Blow-up and blow-down}

Suppose that $X$ is a closed symplectic manifold of dimension $2n$
and $Y\subset X$ is a submanifold of $X$ of codimension $2k$.

 As we are going to apply
the degeneration formula, it is convenient to phrase it in terms of
symplectic cut \cite{Le}, which we now review.

\subsubsection{Symplectic cut and normal connected sum}
 Suppose that $X_0\subset X$ is an open codimension
zero submanifold with a Hamiltonian $S^1-$action. Let $H:X_0\to
\mathbb R$ be a Hamiltonian function with $0$ as a regular value. If
$H^{-1}(0)$ is a separating hypersurface of $X_0$, then we obtain
two connected manifolds $X_0^{\pm}$ with boundary $\partial
X_0^{\pm}=H^{-1}(0)$, where the $+$ side corresponds to $H<0$.
Suppose further that $S^1$ acts freely on $H^{-1}(0)$. Then the
symplectic reduction $Z=H^{-1}(0)/S^1$ is canonically a symplectic
manifold of dimension $2$ less. Collapsing the $S^1-$action on
$\partial X^{\pm}=H^{-1}(0)$, we obtain closed smooth manifolds
$\overline{X}^{\pm} $ containing respectively real codimension $2$
submanifolds $Z^{\pm}=Z$ with opposite normal bundles. Furthermore
$\overline {X}^{\pm}$ admits a symplectic structure $\overline
\omega^{\pm}$ which agrees with the restriction of $\omega$ away
from $Z$, and whose restriction to $Z^{\pm}$ agrees with the
canonical symplectic structure $\omega_Z$ on $Z$ from symplectic
reduction. The pair of symplectic manifolds $(\overline {X}^{\pm},
\overline \omega^{\pm})$ is called the symplectic cut of $X$ along
$\Upsilon$.

This is neatly shown by considering $X_0\times \mathbb C$ equipped
with appropriate  product symplectic structures and
 the product $S^1-$action on $X_0\times
\mathbb C$ where $S^1$ acts on $\mathbb C$ by complex
multiplication.
 The extended
action is Hamiltonian if we use the standard symplectic structure
$\sqrt {-1}dw\wedge d\bar w$ or its negative on the $\mathbb C$
factor.

The normal connected sum operation (\cite{G}, \cite{MW}), or the
fiber sum operation is the inverse operation of the symplectic cut.
Given two symplectic manifolds containing symplectomorphic
codimension 2 symplectic submanifolds with opposite normal bundles,
the normal connected sum operation produces a new symplectic
manifold by identifying the tubular neighborhoods.

Notice that we can apply the normal connected sum operation to the
pairs
$$(\overline
X^+, \overline \omega^+|Z^+),(\overline X^-, \overline
\omega^-|Z^-)$$ to recover $(X, \omega)$.

\subsubsection{Blow-up and blow-down}\label{topology}
Now we apply the symplectic cut to construct the blow-up along $Y$.

The normal bundle $N_Y$ is a symplectic vector bundle, i.e. a bundle
with fiber $(\mathbb R^{2k}, \omega_{std})$. Picking a compatible
almost complex structure on $N_Y$, we then have an Hermitian bundle.
Let $P$ be the principal $U(k)$ bundle.

Now pick a unitary connection $A$ for $P$, and let $\mathcal
W_A\subset u(k)^*$ be as in Lemma \ref{W}. Let
$D_{\epsilon_0}\subset \mathbb C^k=\mathbb R^{2k}$ be the closed
$\epsilon_0-$ball such that its image under  the moment map lies
inside $\mathcal W_A$.

Apply the universal construction to $P$ and  $D_{\epsilon_0}\subset
\mathbb C^k$, we get a symplectic  form $\omega_{\epsilon_0, A}$ on
the disc bundle $N_Y(\epsilon_0)$ which restricts to $\omega_{std}$
on each fiber, and restricts to $\omega|_Y$ on the zero section.

 By the symplectic neighborhood
theorem, and by possibly taking a smaller $\epsilon_0$, a tubular
neighborhood ${\mathcal N}_{\epsilon_0}(Y)$ of $Y$ in $X$ is
symplectomorphic to the disc bundle $N_Y(\epsilon_0)$ with the
symplectic form $\omega_{\epsilon_0, A}$. Let $\phi:{\mathcal
N}_{\epsilon_0}(Y)\to N_Y(\epsilon_0)$ be such a symplectomorphism.

Consider the Hamiltonian $S^1-$action on $X_0={\mathcal
N}_{\epsilon_0}(Y)$ by complex multiplication. Fix $\epsilon$ with
$0< \epsilon<\epsilon_0$ and  consider the moment map
$$H(u)=|\phi(u)|^2-\epsilon, \quad u\in \mathcal N_Y(\epsilon_0),$$
where $|\phi(u)|$ is the norm of $\phi(u)$ considered as a vector in
 a fiber of the Hermitian bundle $N_Y$. Here $X_0\times \mathbb C$ is just
$${\mathcal
N}_{\epsilon_0}(Y)\oplus {\mathbb C}.$$ Write  the hypersurface
$P=H^{-1}(0)$ in $X$ corresponding to the sphere bundle of $N_Y$
with radius $\epsilon$.

We cut $X$ along $P$ to obtain two closed symplectic manifolds
$\overline X^+$ and $\overline X^-$. Notice that $Y$ is contained in
the $+$ side.

$\overline X^-$ is called the blow-up of $X$ along $Y$. Notice that
the construction depends on $\epsilon$, the connection $A$ and the
symplectomorphism $\phi:{\mathcal N}_{\epsilon_0}(Y)\to
N_Y(\epsilon_0)$. However, as remarked in p. 250 in \cite{MS}, given
two different choices $A, \phi$ and $A', \phi'$, for sufficiently
small $\epsilon$, the resulting symplectic forms are isotopic. We
will often denote the blow-up by $\tilde X$ ignoring various
choices.

Denote the codimension 2 symplectic submanifold $Z$ by $E$. We will
call $E$ the exceptional divisor.

Blowing  down is the inverse operation of blowing up. More
precisely, if $\tilde X$ is the blow up of $X$ along $Y$ with the
exceptional divisor $E$. Then, as remarked,
 the normal connected sum of $\tilde X$ and
$\overline X^+$ along $E$ gives back $X$. This process from $\tilde
X$ to $X$ is called the blow down along $E$.

Now we set up to describe the topology of the blow-up.

As smooth manifolds, $E$ is diffeomorphic to the projectivization of
$N_Y$, and
  $\overline{X}^+ $ is the projectivization of
$N_Y\oplus {\mathbb C}$.

Observe that
$$\tilde X=(X-{\mathcal
N}_{\epsilon}(Y))\cup_{\phi} \overline {N_Y(\epsilon_0)}^-.$$ We can
define a map
\begin{equation}p:\tilde X\to X\label{p}
\end{equation}
 which is identity away from ${\mathcal N}_{\epsilon_0}(Y)$. Such
a map can be constructed by identifying ${\mathcal
N}_{\epsilon_0}(Y)-{\mathcal N}_{\epsilon}(Y)$ with the deleted
neighborhood ${\mathcal N}_{\epsilon_0}(Y)-Y$ using a diffeomorphism
from $(\epsilon, \epsilon_0)$ to $(0, \epsilon_0)$. Such a $p$ is
not unique, but the induced maps  $p_*$ and $p^*$ on homology and
cohomology are the same for different choices.

In particular, if $Y$ is of codimension 2, then $E=Y$. And $\tilde
X$ is diffeomorphic to $X$, although the symplectic structures are
not quite the same.

It is important to observe that the pair $(\tilde X, E)$ is the
common piece of the symplectic cuts of $X$ and $\tilde X$. More
explicitly,
 $X$ degenerates into $(\tilde X, E)$ and
$({\mathbb P}(N_{Y}\oplus {\mathbb C}), E)$, while $\tilde X$
degenerates into $(\tilde X, E)$ and $({\mathbb P}(N_{E|\tilde
X}\oplus {\mathbb C}), E)$.

Finally we mention that it is also shown in \cite{GS} that blowing
up/down can be explicitly realized as a birational cobordism.


\section{Relative GW invariants and the degeneration formula}

Li and Ruan \cite{LR} first introduced the moduli space of relative
stable maps and constructed its virtual fundamental class.
Integrating against the virtual fundamental class, they first
defined the relative Gromov-Witten invariants (see \cite{IP} for a
different version and \cite{Li} for the algebraic treatment). They
are the main tool of the paper. We want to review briefly the
construction.

\subsection{GW-invariants}\label{gw}

Suppose that $(X,\omega)$ is a compact symplectic manifold  and $J$
is a tamed almost complex structure.

\begin{definition} \label{stable map}A stable $J-$holomorphic map is an equivalence class of
pairs $(\Sigma, f)$. Here $\Sigma$ is a connected nodal marked
Riemann surface with arithmetic genus $g$, $k$ smooth marked points
$x_1,  \cdots , x_k$, and $f : \Sigma \longrightarrow X$ is a
continuous map whose restriction to each component of $\Sigma$
(called a component of $f$ in short) is $J$-holomorphic.
Furthermore, it satisfies the stability condition: if $f\mid _{S^2}$
is constant (called a ghost bubble) for some $S^2-$component, then
the  $S^2-$component has at least three special points (marked
points or nodal points). $(\Sigma, f)$, $(\Sigma', f')$ are
equivalent, or $(\Sigma,f) \sim (\Sigma', f')$, if there is a
biholomorphic map $h : \Sigma' \longrightarrow \Sigma$ such that
$f'=f\circ h$.
\end{definition}

An essential feature of Definition \ref{stable map} is that, for a
stable $J-$holomorphic map $(\Sigma, f)$, the automorphism  group
$$
   \mbox{Aut}(\Sigma, f) = \{ h\mid h\circ (\Sigma, f) = (\Sigma, f)\}
$$
is finite. We define the moduli space $\overline{\mathcal
M}^X_A(g,k, J)$ to be the set of equivalence classes of stable
$J-$holomorphic maps such that $[f] = f_*[\Sigma] =A\in H_2(X,{\bf
Z})$. The virtual dimension of $\overline{\mathcal M}^X_A(g,k,J)$ is
computed by index theory,
$$
\mbox{virdim}_{\mathbb R}\overline{\mathcal M}^X_A(g,k,J) = 2c_1(A)
+ 2(n-3)(1-g) +2k,
$$
where $n$ is the complex dimension of $X$.

Unfortunately, $\overline{\mathcal M}^X_A(g,k,J)$ is highly singular
and may have larger dimension than the virtual dimension. There are
several ways to extract invariants (\cite{FO}, \cite{LT}, \cite{S1},
\cite{R1}), the one we use is the  virtual neighborhood method in
\cite{R1}.

First, we drop the $J$-holomorphic condition from the previous
definition and require only each component of $f$ be smooth. We call
the resulting object a stable map or a $C^\infty$-stable map. Denote
the corresponding space of equivalence classes by
$\overline{\mathcal B}^X_A(g,k,J)$. $\overline{\mathcal
B}^X_A(g,k,J)$ is clearly an infinite dimensional space. It has a
natural stratification given by the topological type of $\Sigma$
together with the fundamental classes of the components of  $f$. The
stability condition ensures that $\overline{\mathcal B}^X_A(g,k,J)$
has only finitely many strata such that each stratum is a Frechet
orbifold. Further one can use the pregluing construction to define a
topology on $\overline{\mathcal B}^X_A(g,k,J)$ which is Hausdorff
and makes $\overline{\mathcal M}^X_A(g,k,J)$ a compact subspace (see
\cite{R1}).

We can define another infinite dimensional space $\Omega^{0,1}$
together with a map
$$
  \pi : \Omega ^{0,1} \longrightarrow \overline{\mathcal
B}^X_A(g,k,J)
$$
such that the fiber is $\pi^{-1}(\Sigma,f) = \Omega ^{0,1} (f^*TX)$.
The Cauchy-Riemann operator is now interpreted as a section of $ \pi
: \Omega ^{0,1} \longrightarrow \overline{\mathcal B}^X_A(g,k,J)$,
$$\overline{\partial}_J:\overline{\mathcal B}^X_A(g,k,J)\to \Omega ^{0,1}$$
with $\overline{\partial}_J^{-1}(0)$  nothing but
$\overline{\mathcal M}^X_A(g,k,J)$.

 At each $(\Sigma,f) \in
\overline{\mathcal M}^X_A(g,k,J)$, there is a canonical
decomposition of the tangent space of $\Omega^{0,1}$ into the
horizontal piece and the vertical piece. Furthermore,  by choosing a
compatible Riemannian  metric on $X$ we can linearize
$\overline{\partial}_J$ with respect to deformations of stable maps
and project to the vertical piece to obtain an elliptic complex over
$\overline{\mathcal B}^X_A(g,k,J)$,
\begin{equation}\label{virtual-0}
    L_{\Sigma,f}: \Omega^0(f^*TX) \longrightarrow
    \Omega^{0,1}(f^*TX).
\end{equation}

Since $\overline{\mathcal M}^X_A(g,k,J)$ is compact, we can
construct a global orbifold bundle $\mathcal E$ over
$\overline{\mathcal B}^X_A(g,k,J)$ together with a bundle map
$\eta:\mathcal E\to \Omega^{0,1}$ supported in a neighborhood
$\mathcal U$ of $\overline{\mathcal M}^X_A(g,k,J)$.

Consider  the finite dimensional vector bundle  over ${\mathcal U}$,
$p: {\mathcal E}|_{\mathcal U} \longrightarrow {\mathcal U}$. The
stabilizing equation $\overline{\partial}_J + \eta$ can be
interpreted as a section of the bundle $p^*\Omega^{0,1}\to {\mathcal
E}|_{\mathcal U}$. By construction this section
$$
\overline{\partial}_J + \eta: {\mathcal E}|_{\mathcal U}\to
p^*\Omega^{0,1}$$ is transverse to the zero section of
$p^*\Omega^{0,1}\to {\mathcal E}|_{\mathcal U}$.

The set $U^X_{{\mathcal S}_e} = (\overline{\partial}_J
+\eta)^{-1}(0)$    is called a virtual neighborhood in [R1]. The
heart of \cite{R1} is to show that $U^X_{{\mathcal S}_e} $ has the
structure of a $C^1-$manifold.

Notice that $U^X_{{\mathcal S}_e} \subset {\mathcal E}\mid_{\mathcal
U}$.
 Over  $U^X_{{\mathcal S}_e}$ there is the tautological bundle
 $${\mathcal E}_X =
p^*({\mathcal E}\mid _U)|_{U^X_{{\mathcal S}_e}}. $$ It comes  with
the tautological inclusion map
$$
  S_X: U^X_{{\mathcal S}_e} \longrightarrow {\mathcal E}_X, \quad ((\Sigma', f'), e)\to e,
$$
which can be viewed as a section of ${\mathcal E}_X$. It is easy to
check that
$$S^{-1}_X(0) = \overline{\mathcal M}^X_A(g,k,J).$$
Furthermore, one can show that $S_X$ is a proper section.

There are evaluation maps
$$
  ev_i: \overline{\mathcal B}_A(g,k,J)\longrightarrow X, \quad
  (\Sigma, f) \to f(x_i),
$$
for $1\leq i\leq k$. $ev_i$ induces a natural map from $U_{{\mathcal
S}_e}\longrightarrow X^k$, which can be shown to be smooth.

Let $\Theta$ be the Thom form of the finite dimensional bundle
${\mathcal E}_X\to U^X_{{\mathcal S}_e}$.

\begin{definition} The (primitive and primary) GW invariant is
defined as
$$
  \langle\alpha_1,\cdots,\alpha_k\rangle^X_{g,A} = \int_{U_{{\mathcal
  S}_e}}S_X^*\Theta\wedge\Pi_i ev_i^*\alpha_i,
$$
where  the $\alpha_i$ are classes in  $H^*(X;{\mathbb R})$ and are
called primary insertions. For the genus zero case, we also
sometimes write $\langle\alpha_1,\cdots,\alpha_k\rangle^X_A$ for $
\langle\alpha_1,\cdots,\alpha_k\rangle^X_{0,A}$.
\end{definition}

\begin{definition} For each marked point $x_i$, we define an orbifold complex line
bundle ${\mathcal L}_i$ over $\overline{\mathcal B}^X_A(g,k,J)$
whose fiber is $T_{x_i}^*\Sigma$ at $(\Sigma,f)$. Such a line bundle
can be pulled back to $U^X_{{\mathcal S}_e}$ (still denoted by
${\mathcal L}_i$). Denote $c_1({\mathcal L}_i)$, the first Chern
class of ${\mathcal L}_i$, by $\psi_i$.
\end{definition}

\begin{definition} The descendent GW invariant is
defined as
$$
  \langle\tau_{d_1}\alpha_1,\cdots, \tau_{d_k}\alpha_k\rangle^X_{g,A}  =
 \int_{U^X_{{\mathcal S}_e}}S_X^*\Theta\wedge \Pi_i\psi_i^{d_i}\wedge ev_i^*\alpha_i,
$$
where $\alpha_i\in H^*(X;{\mathbb R})$.
\end{definition}

\begin{remark}
In the stable range $2g +k \geq 3$, one can also define
non-primitive GW invariants (See e.g. \cite{R1}). Recall that there
is a map $\pi: \overline{\mathcal B}^X_A(g,k,J)\rightarrow
\overline{\mathcal M}_{g,k}$ contracting the unstable components of
the source Riemann surface. We can bring  in a class $\kappa$ from
the Deligne-Mumford space via $\pi$ to define  the ancestor GW
invariants
$$
  \langle \kappa\parallel \Pi_i\alpha_i\rangle^X_{g,A} =
 \int_{U^X_{{\mathcal S}_e}}S_X^*\Theta\wedge \pi^*\kappa\wedge \Pi_i ev_i^*\alpha_i.
$$
The particular class we will use is the point class in
$\overline{\mathcal M}_{0,k}$ (see Theorem \ref{projective}).
\end{remark}

All GW invariants are invariants of $(X, \omega)$. In fact they are
invariant under  deformations $\omega_t$ of $\omega$ by Remark
\ref{deformation2}.

\begin{remark}\label{D} For computational purpose we would
use the following variation of the virtual neighborhood construction
in Section 6.
 Suppose $\iota: D\subset X$ is a submanifold. For $\alpha  \in H^*(D;
{\mathbb R})$ we define  $\iota!(\alpha) \in H^*(X;{\mathbb R})$ via
the transfer map $i! = PD_X\circ \iota_* \circ PD_D$. One can
construct GW-invariants with an insertion of the form $i!(\alpha)$
as follows. Apply the virtual neighborhood construction to the
compact subspace
$$\overline{\mathcal M}_A(g,k,J)\cap ev_1^{-1}(D)$$
in $ \overline{\mathcal B}^X_A(g,k,J,D) = ev_1^{-1}(D)$  to obtain a
virtual neighborhood ${\mathcal U}_{{\mathcal S}_e}(D)$ together
with the natural map $ev_D : {\mathcal U}_{{\mathcal S}_e}(D)
\longrightarrow D$. It is easy to show that
\begin{eqnarray*}
  \langle \tau_{d_1}i!(\alpha), \tau_{d_2}\beta_2,\cdots,
  \tau_{d_k}\beta_k\rangle^X_{g,A}
    =  \int_{U_{{\mathcal
 S}_e}(D)}S^*\Theta\wedge ev_D^*\alpha\wedge \prod_{i=2}^k \psi_i^{d_i}ev_i^* \beta_i
  .
\end{eqnarray*}
\end{remark}

\begin{remark} \label{absolute graph}
        For each $ \langle\tau_{d_1}\alpha_1,\cdots, \tau_{d_k}\alpha_k\rangle^X_{g,A}$, we can conveniently associated a simple graph $\Gamma$ of one vertex decorated by
        $(g,A)$ and   a tail for each marked point. We then further decorate
        each tail by $(d_i, \alpha_i)$ and call the  resulting graph $\Gamma(\{(d_i,\alpha_i)\})$ {\em a weighted
        graph}. Using the weighted graph
        notation, we denote the above invariant by $\langle \Gamma(\{(d_i, \alpha_i)\})\rangle^X$.
        We can also consider the
       disjoint union $\Gamma^{\bullet}$ of several such graphs  and use $A_{\Gamma^{\bullet}},
       g_{\Gamma^{\bullet}}$ to denote the total homology class and total arithmetic genus.
       Here the total arithmetic genus is $1+\sum (g_i-1)$.  Then,
       we define $\langle \Gamma^{\bullet}(\{(d_i, \alpha_i)\})\rangle^X$ as the  product
       of GW invariants of the connected components.
       \end{remark}

\subsection{Relative GW-invariants}\label{rgw}

In this section, we will review the relative GW-invariants.
 The readers can
find more details in the reference \cite{LR}.

Let $Z\subset X$ be a real codimension $2$ symplectic submanifold.
Suppose that $J$ is an $\omega-$tamed almost complex structure on
$X$ preserving $TZ$, i.e. making $Z$ an almost complex submanifold.
The relative GW invariants are defined by counting the number of
stable $J-$holomorphic maps intersecting $Z$ at finitely many points
with prescribed tangency. More precisely, fix a $k$-tuple $T_k=(t_1,
\cdots, t_k)$ of positive integers, consider a marked pre-stable
curve
$$
(C,x_1, \cdots, x_m, y_1, \cdots, y_k)
$$
and stable $J-$holomorphic maps $
 f: C\longrightarrow X
$ such that the divisor $f^*Z$ is
$$
f^*Z = \sum_i t_i y_i.
$$
One would like to consider the moduli space of such curves and apply
the virtual neighborhood technique to construct the relative
invariants. But this scheme needs modification as the moduli space
is not compact. It is true that for a sequence of $J-$holomorphic
maps $(\Sigma_n,f_n)$ as above, by possibly passing to a
subsequence, $f_n$ will still converge to a stable $J$-holomorphic
map $(\Sigma,f)$. However the limit $(\Sigma,f)$ may have some
$Z-$components, i.e. components whose images under $f$ lie entirely
in $Z$.

To deal with this problem  the authors in \cite{LR} adopt  the open
cylinder model.  Choose a Hamiltonian $S^1$ function $H$ in a closed
$\epsilon-$symplectic tubular neighborhood $X_0$ of $Z$ as in 2.2.2
with $H(X_0)=[-\epsilon, 0]$ and $Z=H^{-1}(-\epsilon)$. Next we need
 to choose an almost complex structure with nice properties near $Z$.
An almost complex structure $J$ on $X$ is said to be tamed relative
to $Z$  if $J$ is $\omega$-tamed, $S^1-$invariant for some $(X_0,
H)$, and  such that $Z$ is an almost complex submanifold. The set of
such $J$ is nonempty and forms a contractible space. With such a
choice of almost complex structure, $X_0$ can be viewed as a
neighborhood of the zero section of the complex line bundle
$N_{Z|X}$ with the $S^1$ action given by  the complex multiplication
$e^{2\pi i \theta}$. Now we remove $Z$. The end of $X-Z$ is simply
$X_0-Z$. Recall that the punctured disc $D-\{0\}$ is biholomorphic
to the half cylinder $S^1\times [0, \infty)$. Therefore, as an
almost complex manifold, $X_0-Z$ can be viewed as the translation
invariant almost complex half cylinder $P\times  [0, \infty)$ where
$P=H^{-1}(0)$. In this sense, $X-Z$ is viewed as a manifold with
almost complex cylinder end.

Now we consider a holomorphic map in the cylinder model where the
marked points intersecting $Z$ are removed from the domain surface.
Again we can view a punctured neighborhood of each of these marked
points as a half cylinder $S^1 \times [0,\infty)$. With such a $J$,
a $J-$holomorphic map of $X$ intersecting $Z$ at finitely many
points then exactly corresponds to a $J-$holomorphic map to the open
manifold $X-Z$ from a punctured Riemann surface which converges to
(a multiple of) an $S^1-$orbit at a punctured point.

Now we reconsider the convergence of $(\Sigma_n,f_n)$ in the
cylinder model. The creation of a $Z-$component $f_\nu$ corresponds
to disappearance of a part of $im(f_n)$ into the infinity. We can
use translation to rescale back the missing part of $im(f_n)$. In
the limit, we may obtain a stable map $\tilde{f}_\nu$ into $P\times
{\mathbb R}$. When we obtain $X$ from the cylinder model, we need to
collapse the $S^1$-action at the infinity. Therefore, in the limit,
we need to take into account  maps into the closure of $P\times
{\mathbb R}$. Let $Q$ be the projective completion of the normal
bundle $N_{Z|X}$, i.e. $Q = {\mathbb P}(N_{Z|X}\oplus {\mathbb C})$.
Then $Q$ has a zero section $Z_0={\mathbb P}(0\oplus {\mathbb C})$
and an infinity section $Z_\infty={\mathbb P}(N_{Z|X}\oplus 0)$. One
can further show that $\tilde{f}_\nu$ indeed is a stable map into
$Q$ with the stability specified below.

To form a compact moduli space of such maps we thus must allow the
target $X$ to degenerate as well (compare with \cite{Li}). For any
non-negative integer $m$, construct $Q_m$ by gluing together $m$
copies of $Q$, where the infinity section  of the $i^{th}$ component
is glued to the zero section  of the $(i+1)^{th}$ component for
$1\leq i \leq m$. Denote the zero section of the $i^{th}$ component
by $Z_{i, 0}$, and the infinity section by $Z_{i,\infty}$, so  $Sing
Q_m = \cup_{i=1}^{m-1} Z_{i,\infty}$. We will also denote
$Z_{m,\infty}$ by $Z_\infty$ if there is no possible confusion.
Define $X_m$ by gluing $ X $  to $Q_m$ along $Z\subset X$ and
$Z_{1,0}\subset Q_m$. In particular, $X_0=X$ will be referred to as
the root component and the other irreducible components will be
called the bubble components.

$Z\subset X$ can be thought of as the infinity section
$Z_{0,\infty}$ of the $0-$th component (which is $X$), thus $Sing
X_m = \cup_{i=0}^{m-1}Z_{i, \infty}$. Let $\mbox{Aut}_Z Q_m$ be the
group of automorphisms of $Q_m$ preserving $Z_{i,0}$,
$Z_{i,\infty}$, and the morphism to $Z$.  And let $\mbox{Aut}_ZX_m$
be the group of automorphisms of $X_m$ preserving $X$ (and $Z$) and
with restriction to $Q_m$ being contained in $\mbox{Aut}_ZQ_m$ (so
$\mbox{Aut}_ZX_m = \mbox{Aut}_ZQ_m\cong ({\mathbb C}^*)^m$, where
each factor of $({\mathbb C}^*)^m$ dilates the fibers of the $i-$th
${\mathbb P}^1-$bundle). Denote by $\pi_m : X_m\longrightarrow X$
the map which is the identity on the root component $X_0$ and
contracts all the bubble components to $Z=Z_{1,0}$ via the  fiber
bundle projections.

Now consider a nodal curve $C$ mapped into $X_m$ by
$f:C\longrightarrow X_m$ with specified tangency to $Z_{m,\infty}$.
There are two types of marked points:

(i) absolute marked points whose image under $f$ lie outside
$Z_{m,\infty}$ labeled by $x_i$,

(ii) relative marked points which are mapped into $Z_{m,\infty}$ by
$f$ labeled by $y_j$.

A relative $J-$holomorphic map $f:C\longrightarrow X_m$ is said to
be pre-deformable if $f^{-1}(Z_{i-1,\infty}=Z_{i,0})$ consists of a
union of nodes so that for each node $p\in
f^{-1}(Z_{i-1,\infty}=Z_{i,0}), i=1,2, \cdots, m$, the two branches
at the node are mapped to different irreducible components of $X_m$
and the orders of contact to $Z_{i-1,\infty}=Z_{i,0}$ are equal.

An isomorphism of two such $J-$holomorphic maps $f$ and $ f'$ to
$X_m$ consists of a diagram
$$
\begin{array}{ccc}
     (C,x_1, \cdots, x_l, y_1, \cdots, y_k) &
     \stackrel{f}{\longrightarrow} &  X_m\\
     h\downarrow &  & \downarrow t\\
     (C',x'_1, \cdots, x'_l, y'_1, \cdots, y'_k) &
     \stackrel{f'}{\longrightarrow} &  X_m
\end{array}
$$
where $h$ is an isomorphism of marked curves and $t\in
\mbox{Aut}_Z(X_m)$. With the preceding understood, a relative
$J-$holomorphic map to $X_m$ is said to be stable if it has only
finitely many automorphisms.

We introduced the notion of a weighted graph in Remark \ref{absolute
graph}. We need to refine it for relative stable maps to $(X,Z)$. A
(connected) relative graph $\Gamma$ consists of the following data:

(1) a vertex decorated by $A\in H_2(X;\mathbb Z)$ and genus $g$,

(2) a tail for each absolute marked point,

(3) a relative tail for each relative marked point.

\begin{definition} Let $\Gamma$ be a relative graph
with $k$ (ordered) relative tails  and $T_k = (t_1, \cdots, t_k)$, a
$k-$tuple of positive integers forming a
 partition of $A\cdot [Z]$.
 A relative
$J-$holomorphic map of $(X, Z)$ with  type $(\Gamma, T_k)$ consists
of a marked curve $(C, x_1, \cdots, x_l, y_1, \cdots, y_k)$ and a
map $f: C\longrightarrow X_m$ for some non-negative integer $m$ such
that


(i) $C$ is a connected curve (possibly reducible) of arithmetic
genus $g$,
 (ii) the map
$$
   \pi_m\circ f : C \longrightarrow X_m \longrightarrow X
$$
satisfies $(\pi_m\circ f)_*[C] = A$,

(iii) the  $x_i, 1\leq i \leq l$, are the absolute marked points,

(iv) the $y_i, 1\leq i \leq k$, are the relative marked points,

(v) $f^{-1}(Z_{m, \infty})$ consists of precisely the points
$\{y_1,\cdots, y_k\}$ and $f$ has order $t_i$ at each $y_i$.

\end{definition}

Let $\overline{\mathcal M}_{\Gamma, T_k}(X,Z,J)$ be the moduli space
of equivalence classes of pre-deformable relative stable
$J-$holomorphic maps with type $(\Gamma, T_k)$. Notice that for an
element $f:C\to X_m$ in $\overline{\mathcal M}_{\Gamma, T_k}(X,Z,J)$
the intersection pattern with $SingX_m$ is only constrained by the
genus condition and the pre-deformability condition.

Consider the configuration space $\overline{\mathcal B}_{\Gamma,
T_k}(X,Z,J)$ of equivalence classes of smooth pre-deformable
relative stable maps to $X_m$ for all $m\geq 0$. Here, for each $m$,
we still take the equivalence class under Aut$_ZX_m$. In particular,
the subgroup of Aut$_ZX_m$ fixing such a map is required to be
finite. The maps are required to intersect the $Z_{i,\infty}$ only
at finitely many points in the domain curve. Further, at these
points, the map is required to have a holomorphic leading term  in
the normal Taylor expansion for any local chart of $X$ taking $D$ to
a coordinate hyperplane and being holomorphic in the normal
direction along $D$. Thus the notion of contact order still makes
sense, and we can still impose the
 pre-deformability condition and contact order condition at the $y_i$
 being  governed by $T_k$.

With the preceding understood, by choosing a unitary connection on
the normal complex line bundles of the $Z_{i,\infty}$,  we can
define the analog of (\ref{virtual-0}),
$$
  L^{X,Z}_{\Sigma, f} : \Omega^0_r \longrightarrow \Omega^{0,1}_r,
$$
taking into account the pre-deformability condition and the contact
condition $T_k$ along $Z_{m,\infty}$. Now we can apply the virtual
neighborhood technique to construct $U^{X,Z}_{{\mathcal S}_e}$,
$\mathcal E_{X,Z}$, $S_{X,Z}$ as in section \ref{gw}.

In addition to the  evaluation maps on $\overline{\mathcal
B}_{\Gamma, T_k}(X,Z,J)$,
$$\begin{array}{lllll}
 ev^X_i: &\overline{\mathcal B}_{\Gamma, T_k}(X,Z,J)& \longrightarrow &X, &\quad 1\leq i\leq l,\\
&f&\longrightarrow &\pi_m\circ f(x_i),&
\end{array}
$$
there are also the evaluations maps
$$\begin{array}{lllll}
  ev^Z_j: &\overline{\mathcal B}_{\Gamma, T_k}(X,Z,J) &\longrightarrow& Z, &\quad  1\leq j\leq
  k,\\
  &f&\longrightarrow & f(y_j),&
\end{array}
$$
where $Z=Z_{m,\infty}$ if the target of $f$ is $X_m$.

\begin{definition} Let $\alpha_i\in H^*(X;{\mathbb R})$,$1\leq i\leq l$, $\beta_j\in
H^*(Z;{\mathbb R})$, $1\leq j\leq k$. Define the  relative GW
invariant
$$
   \langle\Pi_i\tau_{d_i}\alpha_i\mid \Pi_j\beta_j\rangle^{X,Z}_{\Gamma, T_k} =
   \frac{1}{|\hbox{Aut}(T_k)|}\int_{U^{X,Z}_{{\mathcal
   S}_e}}S^*_{X,Z}\Theta\wedge\Pi_i\psi_i^{d_i}\wedge(ev^X_i)^*\alpha_i\wedge\Pi_j(ev^Z_j)^*\beta_j,
$$
where $\Theta$ is the Thom class of the bundle ${\mathcal E}_{X,Z}$
and Aut$(T_k)$ is the symmetry group of the partition $T_k$. Denote
by ${\mathcal T}_k = \{ (t_j, \beta_j)\mid j=1, \cdots, k\}$ the
weighted partition of $A\cdot [Z]$. If the vertex of $\Gamma$ is
decorated by $(g, A)$,  we will sometimes write
$$
   \langle\Pi_i\tau_{d_i}\alpha_i\mid {\mathcal
  T}_k\rangle^{X,Z}_{g,A}
  $$
  for
  $\langle\Pi_i\tau_{d_i}\alpha_i\mid
\Pi_j\beta_j\rangle^{X,Z}_{\Gamma, T_k} $.
\end{definition}

\begin{remark} In \cite{LR} only invariants  without
descendent classes were considered. But it is straightforward  to
extend the definition of \cite{LR} to include absolute descendent
classes.

\end{remark}

      We can  decorate
      the tail of a relative graph $\Gamma$ by $(d_i, \alpha_i)$ as in the absolute case. We can
      further  decorate the relative tails by the weighted partition ${\mathcal T}_k$.
      Denote  the resulting weighted
      relative graph by $\Gamma\{(d_i, \alpha_i)\}| {\mathcal T}_k$. In \cite{LR} the source curve is required
to be connected. We will also need to use a disconnected version.
For a
 disjoint union $\Gamma^{\bullet}$  of weighted relative
      graphs and a corresponding disjoint union of partitions, still denoted by $T_k$,
      we  use $\langle \Gamma^{\bullet}\{(d_i, \alpha_i)\} | {\mathcal
      T}_k
      \rangle^{X,Z}$ to denote the corresponding relative invariants with a disconnected domain,
      which is simply the product of the connected relative
      invariants. Notice that although we use $\bullet$ in our
      notation following \cite{MP}, our disconnected invariants are different.
      The disconnected invariants there depend only on the genus,
      while ours depend on the finer graph data.

\subsection{Degeneration formula}\label{df}
Now we describe the degeneration formula of GW-invariants under
symplectic cutting.

As an operation on topological spaces,  the symplectic cut is
essentially collapsing the circle orbits in the hypersurface
$H^{-1}(0)$ to points in $Z$. Thus we have a continuous map
$$\pi:X\to \overline{X}^+\cup_Z\overline{X}^-.$$
As for the symplectic forms, we have $\omega^+\mid_Z =
\omega^-\mid_Z$. Hence, the pair $(\omega^+, \omega^-)$ defines a
cohomology class of $\overline{X}^+\cup_Z\overline{X}^-$, denoted by
$[\omega^+\cup_Z\omega^-]$. It is easy to observe that
\begin{equation}\label{cohomology relation}
   \pi^* ([\omega^+\cup_Z\omega^-]) = [\omega].
\end{equation}
 Let $B\in
H_2(X;{\mathbb Z})$ be in the kernel of
$$
  \pi_* : H_2(X;{\mathbb Z})
\longrightarrow H_2(\overline{X}^+\cup_Z\overline{X}^-; {\mathbb
Z}). $$
 By (\ref{cohomology relation}) we have $\omega(B) =0$.
Such a class is called a vanishing cycle. The geometric
description of the vanishing cycle is rim tori. For each simple
closed curve $\gamma$ in $Z$, $\pi^{-1}(\gamma)$ is a torus in
$H^{-1}(0)$, i. e. rim tori. It is easy to see that each vanishing
cycle can be represented by a rim tori. In particular, in the case
of blowup along a complex codimension two submanifold $S$,
$H^{-1}(0)$ is the sphere bundle of the normal bundle to $S$ in
$X$. Since the fiber is simply connected, there are no rim tori.
Therefore, there are no vanishing cycle.
 For $A\in H_2(X; {\mathbb Z})$ define $[A] = A + \mbox{Ker}
(\pi_*)$ and
\begin{equation}\label{vanishing cycle}
 \langle\Pi_i\tau_{d_i}\alpha_i\rangle^X_{g,[A]}
= \sum_{B\in[A]}\langle\Pi_i\tau_{d_i}\alpha_i\rangle^X_{g,B} .
\end{equation}
Notice that  $\omega$ has constant pairing with any element in
$[A]$.
 It follows from  the Gromov
compactness theorem that there are only finitely many such elements
in $[A]$ represented by $J$-holomorphic stable maps. Therefore, the
summation in (\ref{vanishing cycle}) is finite.

The degeneration formula expresses
$\langle\Pi_i\tau_{d_i}\alpha_i\rangle^X_{g,[A]} $ in terms of
relative invariants of $(\overline{X}^+, Z)$ and $(\overline{X}^-,
Z)$ possibly with disconnected domains.

To begin with, we need to assume that the cohomology class
$\alpha_i$ is of the form
\begin{equation}\label{decomposition} \alpha_i =
\pi^*(\alpha_i^+\cup_Z\alpha_i^-).\end{equation}
 Here $\alpha_i^\pm \in H^*(\overline{X}^\pm; {\mathbb
R})$ are classes with  $\alpha_i^+\mid_Z = \alpha_i^-\mid_Z$ so that
 they give rise to a class $\alpha_i^+\cup_Z\alpha_i^-\in
H^*(\overline{X}^+\cup_Z\overline{X}^-; {\mathbb R})$.

Next, we proceed to write down the degeneration formula.
    We first specify the relevant
topological type of a marked Riemann surface mapped into
$\overline{X}^+\cup_Z\overline{X}^-$ with the following properties:
\begin{enumerate}
\item[(i)] Each connected component is mapped  either into
$\overline{X}^+$ or $\overline{X}^-$ and carries a respective degree
2 homology class;
\item[(ii)] The marked points are not mapped to $Z$;
\item[(iii)] Each point in the domain mapped to $Z$ carries a positive
integer (representing the order of tangency).
\end{enumerate}

By abusing language we call the above data a $(\overline{X}^+,
\overline{X}^-)-$graph.
 Such a graph gives rise to  two relative
graphs of $(\overline{X}^+, Z)$ and $\overline{X}^-, Z)$, each
possibly being disconnected. We denote them by
 $\Gamma^{\bullet}_+$ and $
\Gamma^{\bullet}_-$ respectively. From (iii) we also get  two
partitions $T_+ $ and $T_-$.  We call   a $(\overline{X}^+,
\overline{X}^-)-$graph a degenerate $(g, A, l)-$graph if the
resulting pairs $(\Gamma^{\bullet}_+, T_+)$ and $
(\Gamma^{\bullet}_-, T_-)$ satisfy the following constraints: the
total number of marked points is $l$, the relative tails are the
same, i.e. $ T_+= T_-$, and the identification of relative tails
produces a connected graph of $X$ with  total homology class
$\pi_*[A]$ and arithmetic genus $g$.

Let $\{\beta_a\}$ be a self-dual basis of $H^*(Z;{\mathbb R})$ and
$\eta^{ab} = \int_Z\beta_a\cup \beta_b$. Given $g, A$ and $ l$,
consider a degenerate  $(g, A, l)-$graph.
 Let $T_k=T_+=T_-$ and ${\mathcal T}_k$ be
a weighted partition $\{t_j, \beta_{a_j}\}$. Let $\breve {\mathcal
T}_k =\{t_j, \beta_{a_{j'}}\}$ be the dual weighted partition.

The degeneration formula for
$\langle\Pi_i\tau_{d_i}\alpha_i\rangle^X_{g,[A]}$ then reads as
follows,
$$
\begin{array}{ll}
 & \langle\Pi_i\tau_{d_i}\alpha_i\rangle^X_{g,[A]}
\cr =&\sum
   \langle \Gamma^{\bullet} \{(d_i, \alpha_i^+)\}|   {\mathcal T}_k\rangle^{\overline{X}^+,Z}
   \Delta({\mathcal T}_k) \langle \Gamma^{\bullet}\{(d_i, \alpha_i^-)\}|
    \breve {\mathcal T}_k\rangle^{\overline{X}^-,Z},\cr
\end{array}
$$
where the summation is taken over all degenerate $(g,A, l)-$graphs,
and
$$\Delta({\mathcal T}_k )=\prod_j t_j Aut(T_k).$$

\section{Uniruled manifolds}

In this section all GW invariants are of genus $0$ and we will omit
the subscript for the genus. This convention will be used in section
6 as well.

\subsection{Uniruledness in algebraic geometry}
Let us first recall the notion of uniruledness in algebraic
geometry.

\begin{definition} A smooth projective variety $X$ (over $\mathbb C$) is called (projectively) uniruled
if for every $x\in X$ there is a morphism $f:{\mathbb P}^1\to X$
satisfying $x\in f(\mathbb P^1)$, i.e. $X$ is covered by rational
curves.

\end{definition}

Rational curves on uniruled projective varieties have the following
nice property \cite{KMM,K1}: for a very general point $x$, if
$g:\mathbb P^1\to X$ is a morphism such that $g_*[\mathbb P^1]=A$
and $g(y_0)=x$, then \begin{equation}\label{convex} H^1(\mathbb P^1,
g^*TX)=0.
\end{equation}

The characterization of uniruled varieties is very important in
classification theory. Here we only review the following beautiful
result which is due to Koll\'ar and Ruan. We sketch a proof as it
may not be so well known.

\begin{theorem} \label{projective}(\cite{K}, \cite{R1}): A projective manifold $(X, J, \omega)$ is uniruled if and
only if there exist a homology class $A\in H_2(X;\mathbb Z)$ and
cohomology classes $\alpha_2,  \cdots , \alpha_{p+1}\in
H^*(X;\mathbb R)$ such that
$$
       \langle [pt]\parallel [pt], \alpha_2,  \cdots , \alpha_{p+1}\rangle^X_A\not= 0,
$$
where the first $[pt]$ represents the Poincar\'e dual of the point
class of $\overline{\mathcal M}_{0,k}$.
\end{theorem}

\begin{proof} Suppose there is a nonzero invariant $\langle [pt]\parallel [pt], \alpha_2,  \cdots , \alpha_{p+1}\rangle^X_A$.
Then for the given (integrable) complex structure $J$, through any
point $x$, the corresponding moduli space of stable $A-$rational
curves cannot be empty. Otherwise the invariant is zero. Since each
domain component of a stable rational curve is $\mathbb P^1$, there
is a morphism $g:\mathbb P^1\to X$ such that $x\in g(\mathbb P^1)$.
So $X$ is uniruled.

Suppose $X$ is uniruled. Fix a sufficiently general point $x$ and a
very ample divisor $\textrm H$. Let $g:\mathbb P^1\to X$ be a
morphism such that $x\in g(\mathbb P^1)$. Such a $g$ exists as $X$
is assumed to be uniruled. Furthermore assume that the class
$A=g_*[\mathbb P^1]$ has the property that the pairing $\textrm
H(A)$ with the ample divisor $\textrm H$ is minimal among all such
classes. Then every rational curve through $x$ is irreducible, i.e.
any stable $A-$rational curve through $x$ is of the form $f:\mathbb
P^1\to X$. Therefore the moduli space of $A-$rational curves through
$x$ is the same as  the compactified moduli space of stable
$A-$rational curves through $x$. We denote this compact moduli space
simply by
 ${\mathcal M}_x$.
Moreover if we invoke the property (\ref{convex}) we  conclude that
${\mathcal M}_x$ is smooth of expected dimension.

Let ${\mathcal E}=f^*TX$. Then ${\mathcal E}$ is a convex
holomorphic bundle over $\mathbb P^1$, and, as a complex bundle, it
is independent of $f$ since $f_*[\mathbb P^1]=A$.
 For any $f\in {\mathcal M}_x$, the tangent space
 $T_f{\mathcal M}_x$ is identified with $\{v\in H^0({\mathcal E})|v(x)=0\}$.
Observe that  there are finitely many points $y_2,  \cdots ,
y_{p+1}\in \mathbb P^1$ such that  for any  holomorphic section
$v\in H^0({\mathcal E})$, $v=0$ if and only if $v(y_i)=0$ for each
$2\leq i\leq p+1$.

Consider the holomorphic evaluation map
$$\Xi_p:{\mathcal M}_x\to X^{p}, \quad f\to (f(y_2),  \cdots , f(y_{p+1})).$$
Its differential is simply given by
$$\delta(\Xi_{p})_f(v)=(v(y_2),  \cdots , v(y_{p+1})).$$
Therefore, by our choice of the $y_i$, $\Xi_p$ is a holomorphic
immersion at
$$(g, y_2,\cdots, y_{p+1}).$$ Now $\Xi_{p}({\mathcal
M}_x)\subset X^{p}$ is a compact complex subvariety of the same
dimension as that of ${\mathcal M}_x$. In particular, it represents
a nonzero homology class $[\Xi_{p}({\mathcal M}_x)]\in
H_*(X^{p};\mathbb Z)$. Furthermore,
$$(\Xi_{p})_*([{\mathcal M}_x])=\lambda [\Xi_{p}({\mathcal M}_x)]$$
for some $\lambda>0$. There are cohomology classes $\alpha_2,
\cdots , \alpha_{p+1}$ coming from very ample divisors such that
$$(\prod_{i=2}^{p+1}\alpha_i )([\Xi_{p}({\mathcal M}_x)])\ne 0.$$
Notice that the points $y_i, 2\leq i\leq p+1$ are fixed, so the
invariant
$$\begin{array}{ll} &\langle [pt]\parallel [pt], \alpha_2,  \cdots , \alpha_{p+1}\rangle^X_A\cr
= &\int_{{\mathcal M}_x} \pi^*[\overline {\mathcal M}_{0,p}]
\Xi_{p}^* (\Pi_{i=2}^{p+1}\alpha_i)=\prod_{i=2}^{p+1}\alpha_i
([\Xi_{p}({\mathcal M}_x)])\ne 0.\cr \end{array}$$

\end{proof}

Observe that we can choose some special $\alpha_i$ in the following
fashion. The simplest is $\Omega^l$, where $\Omega$ is a K\"ahler
form of $X^p$ and $l$ is the dimension of $\Xi_p({\mathcal M}_x)$.
We can choose $\Omega$ as the product of a K\"ahler form $\omega$ on
$X$. It follows that each $\alpha_i$ can be chosen to be some power
of $\omega$. Finally, the nonzero invariant of Theorem
\ref{projective} can be decomposed as the sum of products of
invariants with 3 insertions by associativity.  It implies that (see
also Proposition 7.3 in \cite{Lu2})

\begin{cor} \label{projective2}A projective manifold is projectively uniruled if and only if
$$ \langle [pt], \omega^p, \alpha\rangle^X_A\ne 0$$ for some $A\ne 0, p, \alpha$, where $\omega$
is a K\"ahler form.

\end{cor}

\subsection{Symplectic uniruledness}

Let $X$ be a closed symplectic manifold. Let $E$ be a smooth
symplectic divisor, possibly empty.

\begin{definition} Let $A\in H_2(X;{\mathbb Z})$ be a nonzero  class.  $A$ is said to be a uniruled class if
there is a nonzero GW invariant
\begin{equation}\label{invariant}\langle [pt], \alpha_2, \cdots ,
\alpha_k\rangle^X_A,
\end{equation}
where $\alpha_i\in H^*(X;\mathbb R)$.
 $A$ is said to be a uniruled
class relative to a divisor $E$ if there is a nonzero relative
invariant $\langle[pt], \alpha_2, \cdots , \alpha_k\mid \beta_1,
\cdots , \beta_l\rangle_A^{X,E}$  where $\alpha_i\in H^*(X;\mathbb
R)$ and $\beta_j\in H^*(E;\mathbb R)$.
\end{definition}

\begin{definition}
$X$ is said to be (symplectically) uniruled if there is a uniruled
class.
 $(X, E)$ is said to be uniruled  if there is
 a uniruled class relative to $E$.

\end{definition}

\begin{remark}
It is easy to see that  we could well use the more general
disconnected GW invariants  to define this concept. This flexibility
is important for the proof of the birational cobordism invariance.
\end{remark}

This notion has been studied in the symplectic context by G. Lu (see
[Lu1-3].  Notice that, by \cite{Ltj}, it is not meaningful to define
this notion by requiring that there is a symplectic sphere in a
fixed class through every point, otherwise every simply connected
manifold would be uniruled.

\begin{remark} \label{strong}According to Corollary \ref{projective2} a projectively
uniruled manifold is symplectically uniruled, in fact strongly
symplectically uniruled. Here $X$ is said to be strongly uniruled if
there is a nonzero invariant of the form (\ref{invariant}) with
$k=3$.
\end{remark}

In dimension 4 it follows from \cite{Mc}, \cite{LL1}, \cite{LL2},
\cite{LM} that the converse is essentially true. While in higher
dimensions it follows from \cite{G} (see also \cite{Lu3}) that there
are uniruled symplectic manifolds which are not projective, and it
follows from \cite{R3}   that there could be infinitely many
distinct uniruled symplectic structures on a given smooth manifold.


\subsection{Minimal descendent invariants}

In this paper an absolute descendent invariant is called strict if
one of the insertions is of the form $\tau_k(\gamma)$ with $k\geq
1$. In this section we want to replace a GW-invariant in the
definition of uniruledness with only primary insertions and a
connected domain by a descendent GW-invariant with a possibly
disconnected domain. But the point insertion will be kept to be
$[pt]$, rather than a descendent one $\tau_d([pt])$ with $d>0$.

This additional  flexibility is very important in our proof of the
birational cobordism invariance of uniruledness. The key ingredient
is the fact that, as in algebraic geometry (see e.g. \cite{Ge}),
 a (absolute) descendent class $\psi$ in the moduli
space of genus zero stable maps can always be expressed as a sum of
boundary classes. This is true because, on one hand, such a relation
is well known on the Deligne-Mumford moduli space of genus zero
stable curves, and on the other hand, the difference of  the
descendent class from the pull-back of the corresponding descendent
class on the genus zero Deligne-Mumford moduli space (the so called
ancestor class) is a boundary class as well.

The following well-known fact in the case $A=0$ will be often used.

\begin{lemma} \label{0class} Let $A$ be the zero class. Then

a) A nonzero GW   invariant with only primary insertions must have
exactly 3 insertions. In particular, if one insertion is the point
class, then the invariant is essentially of the form $\langle [pt],
[X], [X]\rangle^X_0=1$ (see e.g. p.230 of \cite{MS2}).

b) Any nonzero strictly descendent GW invariant must also have 3
primary insertions with
the total degree  equal to the dimension of $X$. In particular,
there are at least 4 marked points.
\end{lemma}

In the next lemma we assume all the invariants are connected.

\begin{lemma} \label{3point} Let $\gamma, \alpha$ be cohomology
classes of $X$ and $0\ne A\in H_2(X;\mathbb Z)$.
\begin{itemize}
\item
 If $ \langle [pt], \tau_i(\gamma), \tau_j(\alpha)\rangle^X_A\ne 0$ for
some $i\geq 1, j\geq 1$, then there is a homology class $A'\ne 0$
and a cohomology class $\beta$ such that
$$
 \langle [pt],\tau_i(\gamma), \beta\rangle^X_{A'}\ne 0.
$$
\item
 If $\langle [pt], \tau_i(\gamma), \alpha\rangle^X_A\ne 0$ for some $i\geq 1$, then either $X$
is uniruled or
$$
\langle [pt], \tau_{i-1}(\gamma)\rangle^X_A\ne 0.
$$
\item If $
\langle [pt], \tau_{i}(\gamma)\rangle^X_A\ne 0 $ for some $i\geq 1$
then $X$ is uniruled.
\end{itemize}
\end{lemma}
\begin{proof}
 First of all, recall that when there are $k\geq 3$ marked
 points $x_1, \cdots, x_k$,  $\psi_2$ can be expressed as a sum of boundary divisors (see the
 comments in the beginning of the subsection):
\begin{equation} \label{relation}
\psi_2=\sum _{A=A_1+A_2}D_{(1,3), A_1\mid (2), A_2}.
\end{equation}
Here
 $D_{(1,3), A_1\mid (2), A_2}$ denotes the divisor
consisting of all stable maps with at least two components, one in
class $A_1$ containing $x_1$ and  $x_3$, and the other in class
$A_2$ containing $x_2$. (\ref{relation}) is understood as an
identity in the cohomology of a virtual neighborhood of
$\overline{\mathcal M}_A^X(0, k, J)$.

  When $k=3$ it follows from the relation (\ref{relation}) and the splitting axiom that
\begin{eqnarray*}
& &\langle [pt], \tau_i(\gamma), \tau_j(\alpha)\rangle^X_A\\
&=&\sum_{A=A_1+A_2}\sum_{\mu,\nu}  \langle [pt], \tau_j(\alpha),
e_{\mu}\rangle^X_{A_1}g^{\mu\nu} \langle e^{\nu},
\tau_{i-1}(\gamma)\rangle^X_{A_2}.
\end{eqnarray*}

Fix $i\geq 1$.  Let $A'\ne 0$ be a class with $\langle
[pt],\tau_i(\gamma), \tau_j(\alpha)\rangle_{A'}^X\not= 0$ and the
smallest pairing with $\omega$ among all strictly descendent
invariants with  a point insertion and 2 other insertions. Then we
have either
 $A_1=A', A_2=0$, or $A_1=0, A_2=A'$. Therefore
\begin{eqnarray*}
 & & \langle [pt], \tau_i(\gamma), \tau_j(\alpha)\rangle^X_{A'}\\
&=&\sum_{\mu,\nu} \langle [pt], \tau_j(\alpha),
e_{\mu}\rangle^X_{A'}
g^{\mu\nu} \langle e^{\nu}, \tau_{i-1}(\gamma)\rangle^X_0 \\
& &+\sum_{\mu,\nu} \langle [pt], \tau_j(\alpha), e_{\mu}\rangle^X_0
g^{\mu\nu} \langle e^{\nu}, \tau_{i-1}(\gamma)\rangle^X_{A'} .
\end{eqnarray*}

Next, we want to prove that $j=0$. Suppose first that
 $j\geq 1$.  Then by Lemma \ref{0class} $\langle [pt], \tau_j(\alpha), e_\mu\rangle^X_0 =0$. Moreover
$\langle e^{\nu}, \tau_{i-1}(\gamma)\rangle^X_0$ is always zero.
Thus $ \langle [pt], \tau_i(\gamma),$ $ \tau_j(\alpha)\rangle^X_{A'}
=0$. This contradicts the assumption $ \langle [pt], \tau_i(\gamma),
\tau_j(\alpha)\rangle^X_{A'} \not=0$. So we have $j=0$, i.e.
$$
 \langle [pt], \tau_i(\gamma), \alpha\rangle_{A'}^X \not= 0.
$$
This proves the first part of the lemma.

Suppose now that $j=0$. Then either
$$
 \langle [pt], \alpha, e_{\mu}\rangle^X_{A_1}\ne 0
 $$
for some $A_1\ne 0$, or
$$
 \langle [pt], \alpha, e_{\mu}\rangle
^X_0g^{\mu\nu} \langle e^{\nu}, \tau_{i-1}(\gamma)\rangle^X_A \ne 0.
$$

In the former case, since the insertions $\alpha$ and $e^{\mu}$ are
primary, $A_1$ is a uniruled class.

 In the latter case we must have $\alpha=e^{\mu}=[X]$ and $e^{\nu}=pt$.
 Hence we have  $ \langle [pt] , \tau_{i-1}(\gamma)\rangle^X_A\ne 0$.
 Notice that $A$ is actually a uniruled class if $i-1=0$.

  Finally we deal with the last case, a descendent invariant with 2 insertions.
   As $A\ne 0$ there is a class  $\alpha$
 in $H^2(X;{\mathbb R})$ with $\alpha(A)\ne 0$. By the divisor
axiom for a descendent invariant, if $i\geq 1$,
\begin{equation*}
 \alpha(A) \langle [pt], \tau_i(\gamma)\rangle^X_A
= \langle [pt], \tau_i(\gamma), \alpha\rangle^X_A-
 \langle [pt],
\tau_{i-1}(\gamma\cup \alpha)\rangle^X_A .
\end{equation*}
We have either
$$
 \langle [pt], \tau_i(\gamma), \alpha\rangle^X_A\ne 0,
$$
or
$$
 \langle [pt], \tau_{i-1}(\gamma\cup \alpha)\rangle^X_A\ne 0.
$$
In the former case, notice that $\alpha$ is non-descendent, so we
conclude by the 2nd part that, either $X$ is uniruled or $ \langle
[pt], \tau_{i-1}(\gamma)\rangle^X_A\ne 0$. Therefore we find that
either $X$ is uniruled or there is a nonzero invariant with 2
insertions and less descendent power.
 We can repeat this process to show that either $X$ is uniruled, or
 finally $ \langle [pt],\tau_0(\gamma\cup \alpha^l)\rangle^X_A\ne 0$ for some $0\leq l\leq
 i$. But $X$ is obviously uniruled in the last case as well.

 \end{proof}

The main result in this subsection is

\begin{theorem} \label{descendent} A symplectic manifold $X$ is uniruled
if and only if there is a nonzero, possibly disconnected genus zero
descendent GW invariant
\begin{equation}\label{invariant}\langle [pt], \tau_{j_2}(\alpha_2), \cdots ,
\tau_{j_k}(\alpha_k)\rangle^X_A
\end{equation}
such that  the component with the $[pt]$ insertion  has nonzero
curve class.
\end{theorem}

\begin{proof}
First of all notice that an invariant with disconnected domain is
the product of the invariants of the connected components. Thus the
only if part is clear. And to prove the if part we can assume that
the nonzero descendent invariant is actually connected
  and $A\ne 0$.

The case of $k\leq 3$ has been dealt with in Lemma \ref{3point}. We
now examine descendent invariants with $k\geq 4$ insertions.
 Let $A$ be a class with a nonzero invariant
$\langle [pt], \tau_{j_2}(\alpha_2), \cdots ,
\tau_{j_k}(\alpha_k)\rangle^X_A $ with $j_i>0$ for some $i$ and
$k\geq 4$. Assume that $A$ is such a class with the smallest
$\omega(A)$. Further assume that $k$ is the smallest among such an
$A$.

Since $k\geq 4$, we can apply the boundary relation (\ref{relation})
to obtain
\begin{equation}\label{relation-1}
\begin{array}{ll}
& \langle [pt], \tau_{j_2}(\alpha_2), \cdots ,
\tau_{j_k}(\alpha_k)\rangle^X_A  \\
 =&\sum\sum_{A=A_1+A_2}
\sum_{\mu,\nu}
 \langle [pt], \tau_{j_3}(\alpha_3),\tau_{j_{i_1}}(\alpha_{i_1}),  \cdots ,\tau_{j_{i_{k_1}}}(\alpha_{i_{k_1}}),
 e_{\mu}\rangle^X_{A_1}\\
 &g^{\mu\nu}
 \langle e^{\nu},
 \tau_{j_2}(\alpha_2),\tau_{j_{i_{k_1+1}}}(\alpha_{i_{k_1+1}}),
 \cdots ,\tau_{j_{i_{k-2}}}(\alpha_{i_{k-2}})\rangle^X_{A_2},
\end{array}
\end{equation}
where the leftmost   sum is over all partitions of   $\{4, \cdots,
k\}=\{i_1, \cdots, i_{k_1}\}\cup \{i_{k_1+1},\cdots, i_k\}$. By the
minimality of $\omega(A)$ we have either $A_1=A$ or $A_1=0$.

Suppose the sum  with $A_1=0$ in (\ref{relation-1}) is nonzero.
Consider a nonzero product in this sum. For the first factor we must
have $k_1\geq 1, \alpha_3=e^{\mu}=[X]$ by Lemma \ref{0class}.
Therefore $e_{\nu}=[pt]$ and $ k-2-k_1\leq k-3$. Hence the second
factor
$$
\langle [pt],
 \tau_{j_2}(\alpha_2),\tau_{j_{i_{k_1+1}}}(\alpha_{i_{k_1+1}}), \cdots ,\tau_{j_{i_{k-2}}}
 (\alpha_{i_{k-2}})\rangle^X_{A_2=A}
 $$
 has at most $k-1$ many insertions. By our assumption of the minimality of $k$ among all strictly
descendent invariants for the class $A$,  this is a nonzero
non-descendent invariant, which shows that $X$ is uniruled.

The case where the sum with $A_2=0$ in (\ref{relation-1}) is nonzero
is similar. Consider again a nonzero product in this sum. We must
have $k_1+1\leq k$, or equivalently, $k_1\leq k-1$,  by Lemma
\ref{0class}. By the minimality of $k$ we claim as above that there
is a nonzero invariant of the form $ \langle [pt], \alpha_3,
\alpha_{i_1}, \cdots ,\alpha_{i_{k_1}}, e_{\mu}\rangle^X_{A_1=A}$.
Therefore in any case we have that $X$ is uniruled.

\end{proof}

\begin{remark}We can actually replace (\ref{invariant}) in Theorem
\ref{descendent} by a more general descendent invariant of the form
\begin{equation}\label{invariant1}\langle \tau_{j_1}([pt]), \tau_{j_2}(\alpha_2), \cdots ,
\tau_{j_k}(\alpha_k)\rangle^X_A.
\end{equation}
Such a generalization will be useful in our next paper \cite{tLR}.
\end{remark}
\section{Blow-up Correspondence for GW-invariants }

We assume in this section  that  $X$ is a compact symplectic
manifold. Let $S\subset X$ be a symplectic submanifold of $X$ of
codimension $2k$, $\tilde{X}$  the blow up of $X$ along $S$, and $E$
 the exceptional divisor, which is a $\mathbb P^{k-1}-$bundle over $S$. Let
$p:\tilde{X} \to X$ be the map defined in \ref{topology}. In this
section, we will obtain a correspondence between the relative
GW-invariant of $(\tilde X, E)$ and certain absolute GW-invariants
of $X$. Notice that we allow  $E=S$. Thus the correspondence of this
section can be viewed as a generalization of the correspondence of
Maulik-Pandharipande \cite{MP} in the case $k=1$. The main method
used in the proof of our main theorems is the degeneration formula
for GW-invariants reviewed in section \ref{df}.

\subsection{A refined partial ordering}
In this subsection, we want to order the graphs of certain relative
GW invariants of $(\tilde{X}, E)$ following [MP]. The new feature
here is that we refine the order in the case $k>1$.


The partial order is defined in terms of several preliminary partial
orders. We first deal with those involving classes of the $\mathbb
P^{k-1}-$bundle $E$.

Let $\theta_1={\bf 1}, \cdots , \theta_{m_S}$ be a self dual basis
of $H^*(S;{\mathbb R})$ with respect to the intersection pairing,
where $\bf 1$ is the distinguished degree 0 class. $E$ is a $\mathbb
P^{k-1}-$bundle over $S$, so it has a basis of the form
\begin{equation}\label{basisE}
\{\pi_S^*\theta_i\cup [E]^j\}, \quad  1\leq i\leq m_S, 0\leq j\leq
k-1.
\end{equation}
Here $[E]$ is understood to be the first Chern class of the
tautological line bundle over the projective bundle ${\mathbb
P}(N_{E|\tilde X})$. Denote this basis of $E$ by
$\Theta=\{\delta_t\}, 1\leq t\leq km_S$ with $\delta_1=\bf 1$.
Notice that the basis $\Theta$ is still self dual.

 \begin{definition}
A standard (relative)
  weighted partition $\mu$ with length $l(\mu)$ is a partition
weighted by classes of $E$ from $\Theta$, i.e.
$$
   \mu = \{ (\mu_1, \delta_{K_1}), \cdots, (\mu_{\ell(\mu)},
   \delta_{K_{\ell(\mu)}})\},
$$
where $\mu_i$ and $k_i$ are  positive integers with $K_i\leq km_S$.
\end{definition}

\begin{definition}\label{degree} For $\delta=\pi^*\theta \cup [E]^j
\in H^*(E;\mathbb R)$, we define
$$\deg_S(\delta)=\deg(\theta), \quad
\deg_{f}(\delta)=2j.$$
 For a standard weighted partition $\mu$, we define
$$\deg_S(\mu)=\sum_{i=1}^{l(\mu)}\deg_S(\delta_{K_i}).$$
\end{definition}

\begin{definition}
The set of pairs $(m, \delta)$ where $m\in {\mathbb Z}_{>0}$ and
$\delta \in H^*(E;\mathbb R)$ is partially ordered by the following
size relation:
\begin{equation}\label{size}
   (m, \delta) > (m', \delta')
\end{equation}
\begin{itemize}
\item
if $m>m'$, or \item if $m=m'$ and $\deg_S (\delta) > \deg_S
(\delta')$, or \item if  $m=m', \deg_S (\delta) = \deg_S (\delta')$
and $\deg_f (\delta)
> \deg_f (\delta')$.
\end{itemize}
\end{definition}

We may place the pairs of $\mu$ in decreasing order by size, i.e. by
(\ref{size}).

\begin{definition}
A {\it lexicographic } ordering on weighted partitions is then
defined as follows:
$$
   \mu \stackrel{l}{>} \mu'
$$
if, after placing the pairs in $\mu$ and $\mu'$ in decreasing
order by size, the first pair for which $\mu$ and $\mu'$ differ in
size is larger for $\mu$.
\end{definition}

Next we introduce a relevant partial orders on the curve classes  of
$\tilde X$.

\begin{definition}
 A class $A\in H_2(\tilde X;\mathbb Z)$ is called $\omega-$effective if $A$ is represented by a
 pseudo-holomorphic stable map for
 all $\omega-$tamed almost complex structures.
 For effective classes $A$ and $ A'$ in $ H_2(\tilde X;{\mathbb Z})$, we say that
$A'<A$ if  $p_*A-p_*A'\in H_2(X;\mathbb Z)$ has positive pairing
with the symplectic form on $X$.
\end{definition}

  Let $p:\tilde X\to X$ be the continuous blow down map defined as in
(\ref{p}), and
 \begin{equation}\label{basisX}\sigma_1, \cdots ,
\sigma_{m_X} \end{equation}
 be a basis of $H^{*}(X;{\mathbb R})$. Then, as observed in \cite{Mc2}, $p^*:H^*(X;{\mathbb
R})\to H^*(\tilde X; {\mathbb R})$ is an injection, and
\begin{equation}\label{basistildeX}
\begin{array}{lll}
\gamma_{j}&=p^*\sigma_j, &1\leq j\leq m_X,\\
\gamma_{j+m_X}&=\delta_j\cup [E], &1\leq j\leq km_S,\\
\end{array}
\end{equation}
 generate $H^*(\tilde X;{\mathbb
R})$.
 Here we use $[E]$ to denote both the homology class of $E$ and the dual degree 2 cohomology class in
$H^2(\tilde{X};\mathbb R)$. This $[E]$ actually restricts to the
class in $H^2(E;\mathbb R)$ that we called $[E]$ in (\ref{basisE}).
In particular, for $1\leq i\leq km_S$, the Poincar\'e Dual of
$\gamma_{i+m_X}$ in $\tilde X$ is represented by a cycle lying
inside $E$ which is Poincar\'e dual to $\delta_i$ in $E$.

We would order the following type of connected relative invariants.

\begin{definition} \label{standard}
A connected standard  relative
GW invariant of $(\tilde X,E)$ is of the form
$$
 \langle \varpi|\mu\rangle^{\tilde X,E}_{g,A}=\langle \gamma_{L_1}, \cdots ,
\gamma_{L_n}|\mu\rangle^{\tilde X,E}_{g,A},
$$
where $A$ is $\omega-$effective, $\mu$ is a standard weighted
partition with $\sum_j \mu_j\leq A\cdot [E]$, and $L_i\leq m_X$,
i.e. $\gamma_{L_i}=p^*\sigma_{L_i}$.
\end{definition}


 In this section we use
$\Gamma(\varpi)|\mu$, instead of $\Gamma\{(0, \gamma_{L_i})\}|\mu$,
to denote the (connected) relative weighted graph of
 the relative invariant in Definition \ref{standard}.
 We partially order such weighted graphs
 in the following way.

 \begin{definition}\label{order}
$$
 \Gamma(\varpi')|\mu'\quad
\stackrel{\circ}{<}\quad  \Gamma(\varpi)|\mu
$$
if one of the conditions below holds

(1) $A'<A$,

(2) equality in (1) and the arithmetic genus satisfies $g'<g$,

 (3) equality in (1-2) and $\|\varpi'\| < \|\varpi\|$,

(4) equality in (1-3) and $\deg_S (\mu') > \deg_S (\mu)$,

(5) equality in (1-4) and $\mu'\stackrel{l}{>} \mu$,

\noindent where  $\|\varpi\|$ denotes the number of insertions of
$\varpi$.

\end{definition}

All these inequalities are designed so that the dimension of the
moduli space satisfying the larger constraint/condition is larger.
This explains the seemingly strange conditions (4) and (5)
 where the inequalities are reversed.

The relative invariants obtained by taking  disjoint unions
 of the connected invariants in Definition \ref{standard} are
 called the standard relative
 invariants
 of $(\tilde X, E)$ (with a disconnected domain).
 We extend the partial ordering to a
disjoint union $\Gamma^{\bullet}(\varpi)|\mu$ in the obvious way.

\begin{remark}\label{multiplicationorder}
 It is easy but important to observe
that this extended partial order $\stackrel{\circ}{<}$ is preserved
under disjoint union.
\end{remark}

\begin{remark}If we are only interested in genus zero invariants, then
we can replace $g'<g$ in  (2) by the inequality of the number of
connected components, $n'>n$.
\end{remark}

\begin{lemma} Given a standard relative
invariant, there are only finitely many standard
relative invariants lower in the partial ordering. In particular,
there is a minimal standard invariant with nonzero value.
\end{lemma}

\begin{proof} There is a lower bound on $\omega(A)$ among
all effective classes.

\end{proof}

We call a partially order set with the above property {\em lower
bounded}.

\subsection{Trading insertions}
In this subsection, we associate an absolute descendent invariant of
$X$ to each standard relative descendent invariant of $(\tilde X,
E)$. The critical step is to trade a relative insertion with an
absolute descendent insertion of the blown-down manifold. This can
be thought as the relative-absolute correspondence for insertions.

For a relative insertion $(m, \delta)$ with
$\delta=\pi^*\theta_i\cup [E]^j\in H^*(E;{\mathbb R})$, we associate
the absolute descendent insertion
\begin{equation}\label {taud} \tau_{d(m,\delta)}(\tilde \delta)
\end{equation} on $X$ supported on $S$, where
$$
\begin{array}{ll}\tilde \delta&=\theta_i\cup [S], \\
d(m,\delta)&=km-k+\frac{1}{2}\deg_f(\delta)=km-k+j.
\end{array}$$
Here $[S]$ denote the class of a Thom form of the normal bundle to
$S$ in $X$.

Notice that when $k=1$, $S=E$ and $j$ is always zero, so we simply
have
 $$\tilde \delta_i=\delta_i \cup [E], \quad d(m,\delta)=m-1.$$
It is convenient to view $[E]$ as the class of a Thom form of the
normal bundle to $E$ in $\tilde{X}$ supported near
 the symplectic divisor $E$. In terms of homology constraints,
 $\tilde \delta$ and $\delta$ correspond to the same cycle
 lying inside $E$ as previously remarked.

\begin{definition}\label{tilde}
Given a standard (relative) weighted partition $\mu$,
 let
$$d_i(\mu)=d(\mu_i, \delta_{K_i})=k\mu_i-k+\frac{1}{2}\deg_f(\delta_{K_i}),$$
and
\begin{equation}\label{tildemu}
\tilde \mu=\{\tau_{d_1(\mu)}(\tilde \delta_{K_1}), \cdots ,
\tau_{d_{l(\mu)}(\mu)}(\tilde \delta_{K_{l(\mu)}}) \}.
\end{equation}
\end{definition}

\begin{definition}\label{associate}
The absolute descendent invariant associated to a standard relative
invariant
$$ \langle \Gamma^{\bullet}(\varpi)|\mu\rangle^{ \tilde X, E}$$ is
then
$$ \langle \tilde{\Gamma}^{\bullet}(
\varpi, \tilde \mu)\rangle^{X}.$$   Here $\tilde{\Gamma}^{\bullet}$
is an absolute graph, obtained from the relative graph
$\Gamma^{\bullet}$ by  changing the homology class for each vertex
from $A$ to $p_*(A)$. And the insertions
$\gamma_{L_j}=p^*\sigma_{L_j}$ in $\varpi$ are replaced by
$\sigma_{L_j}$.
\end{definition}


We therefore  consider the following absolute descendent invariants
 of $X$.

 \begin{definition}\label{absolute standard}
    An absolute descendent invariant of $X$ is said to be {\em an absolute
descendent invariant relative to $S$} if its descendent insertions
are supported on $S$, i.e. have the form $\tau_d(\tilde \delta)$
where $\tilde \delta=\theta_i\cup [S]$. Such an invariant is called
a colored standard absolute descendent invariant of $X$ relative to
$S$
if
 the insertions are  divided into two collections,
$\varpi$ and $\tilde \mu$, with each insertion in $\varpi$ being of
the form $\sigma_L$,  and  each insertion in $\tilde \mu$  being of
the form $\tau_d \tilde \delta_K$.
\end{definition}

When $k=1$, a colored  absolute descendent invariant of $X$ relative
to $S$ is called admissible if $\sum_j \mu_j = A\cdot [E]$. Notice
that in this case, we simply have $d_i(\mu)=\mu_i-1$, and $S=E$, so
$$\tilde \mu=\{\tau_{\mu_i-1}(\theta_{K_1}\cup [E]), \cdots ,
\tau_{\mu_{l(\mu)}-1}( \theta_{K_{l(\mu)}}\cup [E]) \}.
$$
Therefore $ \langle \tilde \Gamma^{\bullet}(\varpi, \tilde
\mu)\rangle^{ X}$ agrees with the one in \cite{MP}. Thus every
relative invariant gives rise to an admissible absolute invariant
relative to $S$. Notice that, after possibly adding a number of $E$
insertions, every non-descendent absolute invariant is admissible
for at least one separation.

The following observation is crucial.

\begin{lemma} \label{bijection} If $\mu\ne \mu'$, then $\tilde \mu\ne \tilde \mu'$.
Therefore there is a natural bijection between the set of colored
weighted absolute graphs $\tilde{\Gamma}^{\bullet}( \varpi, \tilde
\mu)$ relative to $S$ and the set of weighted relative graphs in
$\tilde X$ relative to the exceptional divisor $E$ in the blow up
$\tilde{X}$ of $X$ along $S$ in the case $k>1$, and it is still
true in the case $k=1$ if we restrict to the admissible ones.
\end{lemma}

\begin{proof} It suffices to show that $d(m,\delta)$ uniquely
determines $m$ and $\deg_f(\delta)$. The point is that $0\leq
\deg_f(\delta)\leq k-1$. Thus $m$ and $\deg_f(\delta)$ are uniquely
determined by the integral division of $d(m,\delta)$ by $k$.

\end{proof}

Hence we can and will  order the set of colored weighted absolute
graphs $\tilde{\Gamma}^{\bullet}( \varpi, \tilde \mu)$, and in the
case $k=1$, the admissible ones,  in the same way as the weighted
relative graphs. \footnote{Notice that it is possible that several
relative invariants correspond to the same absolute invariant,
because some $\sigma_i$ might be expressed as $\theta_j\cup [S]$.}

 Let $I$ be the partially order set of standard weighted relative graphs
 $\Gamma^{\bullet}( \varpi)| \mu$. Consider the
 infinite dimensional vector space ${\mathbb R}^I_{\tilde X, E}$ whose coordinates
 are ordered in the way compatible with the partial order.
 Given a standard weighted graph $\Gamma^{\bullet}(\varpi)|\mu$, we have the  relative
 invariant $\langle \Gamma^{\bullet}(\varpi)|\mu\rangle^{\tilde
 X,E}$. From the numerical values we can form a vector
 $$v_{\tilde X,E}\in {\mathbb R}^I_{\tilde X, E}.$$
 given by the numerical values.
By Lemma \ref{bijection}, $I$ is also the partially ordered set of
colored standard weighted absolute graphs $\tilde{\Gamma}^{\bullet}(
\varpi, \tilde \mu)$ in the case $k>1$ and the admissible ones in
the case $k=1$. Hence we also have the corresponding vector space
${\mathbb R}^I_{X, S}$ and a vector
$$ v_{X,S}\in {\mathbb R}^I_{X, S},$$
 given by the numerical values of  the absolute invariants of $X$ relative to
$S$,
 $\langle\tilde{\Gamma}^{\bullet}(\varpi, \tilde \mu)\rangle^X$.

\subsection {Full correspondence}
In this subsection we  prove the Relative/Absolute correspondence in
the following form.

\begin{theorem} \label{generalization}
There is an invertible lower triangular linear transformation
$$A_S: {\mathbb
R}^I_{\tilde X, E}\rightarrow {\mathbb R}^I_{X, S}$$  such that
(i) the coefficients of $A_S$ are  local in the sense of being
dependent only on $S$ and its normal bundle; (ii)
$$A_S(v_{\tilde X,E})=v_{X,S}.$$
In particular, $v_{\tilde X,E}$ and  $v_{X,S}$ determine each other.

 Moreover, if  $I_{pt}\subset I$ denotes  the subset indexed by the standard relative
invariants of $(\tilde X,E)$  with the first insertion being the
point insertion, then $A_S$ restricts to an invertible lower
triangular transformation from  ${\mathbb R}_{\tilde  X,
E}^{I_{pt}}$ to ${\mathbb R}_{X, S}^{I_{pt}}$.

Finally, if $I_{0, pt}\subset I_{pt}$ denotes the subset of genus
zero invariants, then $A_S$ further restricts to an invertible lower
triangular transformation from  ${\mathbb R}^{I_{0,pt}}_{\tilde X,
E}$ to ${\mathbb R}^{I_{0,pt}}_{X, S}$.
\end{theorem}

\begin{proof} The idea is as follows.
Since the disjoint union preserves the order of graphs, it is enough
to prove such a correspondence  for connected invariants. We hence
consider a connected relative invariant $\langle \Gamma(
\varpi)|\mu\rangle^{\tilde X, E}$ of $(\tilde X, E)$ and the
associated absolute invariant $ \langle \tilde{\Gamma}( \varpi,
\tilde \mu)\rangle^X$ of $X$ relative to $S$. As mentioned in
\ref{topology}, $\tilde X=X^-$, the $-$ piece of the symplectic cut
of $X$ along a normal sphere bundle over $S$.
 Thus we have the degeneration of $X$ into $(\tilde X=\overline X^-, E)$ and
$({\mathbb P}(N_{S|X}\oplus {\mathbb C}),E)=(\overline X^+, E)$.
We apply the degeneration formula to this connected absolute
invariant $\langle \tilde{\Gamma}(\varpi, \tilde \mu)\rangle ^X$
of $X$ distributing all the $\tilde \mu$ insertions to the
${\mathbb P}^k-$bundle side. Then, the degeneration formula can be
immediately interpreted as expressing the absolute invariant
$\langle \tilde{\Gamma}( \varpi, \tilde \mu)\rangle^X$ as a linear
combination of relative invariants of $(\tilde{X}, E)$ with the
coefficients being essentially certain relative invariants of the
projective bundle. With the preferred distribution of insertions,
the original graph $\Gamma(\varpi)|\mu$ turns out to be the
largest weighted relative graph appearing in the linear
combination.

We only prove the assertion for the case with  the first insertion
being the point insertion, i.e. the case of $I_{pt}$. The proof of
the general case is the same (and easier).

The proof will consists of several steps.
 In the first step we set
up how the degeneration formula is applied.

\noindent{\bf Step I--Set up}.
 We begin with a connected standard weighted relative graph
$$\Gamma( [pt], \varpi)|\mu$$ with the vertex decorated by $(g, A)$.
The associated connected colored standard absolute descendent
invariant of $X$ relative to $S$ can be written as $$\langle [pt],
\varpi, \tilde \mu\rangle ^{X}_{g, p_*(A)}.$$

To apply the degeneration formula let us first explicitly make
the preferred distribution of insertions mentioned above.

 The classes $\tilde \delta_{K_i}$ are supported
 on $S$ and so can be represented by forms with support near $S$.
 Hence we just distribute all of them to the $({\mathbb P}(N_{S|X}\oplus {\mathbb C}), E)$ side,
 i.e. we set $\tilde \delta_{K_i}^-=0$ in (\ref{decomposition}).
 Recall $\varpi$ consists of insertions of the form
 $\sigma_{L_j}$ for $L_j\leq m_X$ (cf (\ref{basisX}).
For each such an insertion we set on the $(\tilde X, E)$,
$$\sigma_{L_j}^-=\gamma_{L_j},$$
with an appropriate extension $\sigma_{L_j}^+$ to the positive side.
 In particular, the relative
invariants of $(\tilde X, E)$ appearing in the degeneration formula
are all standard invariants as in Definition \ref{absolute standard}
and so can be ordered.

It is easy to see  that there are no vanishing cycles in this
case, hence $[\pi_*(A)] = \pi_*(A)$. When $k\geq 2$ this is
observed in \cite{LR}. The point is that any degree 2 homology
class $B$ of $X$ can be represented by a surface away from $S$.
Thus $\pi_*(B)$ is represented by the same surface viewed as a
surface in $\tilde X$. We can extend $p:\tilde X\to X$ from
$\tilde X$ to $\tilde X\cup_E \mathbb P^k-$bundle by collapsing
the $\mathbb P^k$ fibers of the $\mathbb P^k-$bundle and still
denote it by $p$. Composing with this $p$, we conclude that
$p_*\pi_*(B)=B$. Hence $\pi_*(B)=0$ if and only if $B=0$. In the
case $k=1$, $\tilde X=X$ and $\pi_*=\phi_*$ where $\phi$ is the
inclusion of $\tilde X$ into its union with the $\mathbb
P^1-$bundle. Thus the composition $p_*\pi_*$ is actually the
identity map on $H_*(X;\mathbb R)$.

Therefore, by the degeneration formula, we have
\begin{equation}\label{versus}
\begin{array}{lll}
& & \langle [pt], \varpi, \tilde \mu\rangle_{g, p_*A}^{  X}\\
& = &\sum  \langle \Gamma^{\bullet} _{-}([pt],
\varpi_1)|\eta\rangle^{\tilde X, E} \Delta(\eta)
 \langle \Gamma^{\bullet} _{+}(\varpi_2, \tilde \mu)|\breve \eta\rangle^{{\mathbb P}(N_{S|X}\oplus {\mathbb C}), E}.
\end{array}
\end{equation}
The sum on the right is over all $(g,p_*(A), \|\varpi\|+\|\tilde
\mu\|+1 )-$graphs, including all distributions of the insertions
$\varpi$ and all standard intermediate weighted partitions $\eta$.
Here $\Delta(\eta)=\prod_r \eta_r\cdot |\mbox{Aut}(\eta)|$, and
$\breve \eta$ is the dual partition of $\eta$.
Since the basis $\{\delta_i\}$ is self-dual, $\breve \eta$ is still
a standard weighted partition.
 The relative GW
invariants on the right are possibly disconnected.

As mentioned our main claim  is that $\Gamma([pt], \varpi)|\mu$
 is the largest weighted
relative graph among (connected or not)weighted relative graphs
appearing in the linear combination (\ref{versus}). Of course we
are only interested in terms with nonzero coefficients. Our
strategy, following [MP], involves finding conditions for which
the relevant relative invariants of $({\mathbb P}(N_{S|X}\oplus
{\mathbb C}), E)$ are nonzero. We use a fibred almost complex
structure $J$ on ${\mathbb P}(N_{S|X}\oplus {\mathbb C})$ to
evaluate such invariants.


Let $f _{-}:C _{-}\rightarrow \tilde X$ and $f _{+}:C
_{+}\rightarrow {\mathbb P}(N_{S|X}\oplus {\mathbb C})$ be elements
of the relative moduli spaces.   Both $C_{-}$ and $C _{+}$ might be
disconnected.

\noindent{\bf Step II--The coefficient of $\Gamma([pt],
\varpi)|\mu$}. In this step we show that the coefficient $C_0$ of
$\Gamma([pt], \varpi)|\mu$ is nonzero.

In this case the splitting  of $p_*A$ is
\begin{equation} \label{splitting} A _{-}=A, \quad
A _{+}=(A\cdot [E])F
\end{equation}
 where $F$ is the fiber class of
${\mathbb P}(N_{S|X}\oplus {\mathbb C})$.

Thus for each connected component on ${\mathbb P}(N_{S|X}\oplus
{\mathbb C})$ side, we are reduced to the following kind of
connected relative invariants of the fiber $({\mathbb P}^k,
{\mathbb P}^{k-1})$,
\begin{equation}\label{P}
\langle \tau_{nd-1-j}[pt]\mid D^j\rangle_{0,dL}^{{\mathbb
P}^k,{\mathbb P}^{k-1}}
\end{equation}
with $j=\deg_f(\delta_{K_i})$, $d$ a positive integer,  and $D$
the hyperplane class of $H^*(\mathbb P^{k-1};\mathbb Z)$.

For the above relative invariants (\ref{P}) of $({\mathbb P}^k,
{\mathbb P}^{k-1})$ the answer is known when $k=1$ (see \cite{OP}).
It follows that  in this case the coefficient $C_0$ is given by (see
\cite{MP})
$$\prod_j{\frac{1}{(\mu_j-1)!}}(A\cdot [E])^{{\bf 1}(\mu)}\ne 0,$$
where ${\bf 1}(\mu)$ is the number of $(1,\delta_1)=(1,{\bf 1})$ in
$\mu$.

The computation of (\ref{P})in the general case is completed in
Theorem \ref{computation}. Consequently,  when $k>1$,  $C_0$ is the
product of rational numbers of the form
$$  \frac{1}{d^{k-j}[(d-1)!]^k}. $$
In particular, $C_0$ is nonzero as well.

 \noindent{\bf Step III--Curve configuration}. In this
step we conclude that  any configuration that might occur with
nonzero coefficient is no bigger  than
 the extremal configuration where $C _{-}$ is a connected genus $g$
 curve,
and $C _{+}$ consists of $l(\eta)$ rational connected components
each having exactly one relative marking (i.e. totally ramified over
the zero section $E_0$).


We first argue the class $A _{+}$ of a largest weighted graph must
be a multiple of the fiber class.
 One splitting of $p_*A$
is $A _{-}=A, A _{+}=(A\cdot [E])F$ where $F$ is the fiber class of
${\mathbb P}(N_{S|X}\oplus {\mathbb C})$. Any other splitting
differs by a class $\beta$ in $H_2(E;\mathbb Z)$. An equivalent
description of the splitting is that
$$p_*A _{-}+ \iota_*t_*A _{+}=p_*A,$$ where $\iota: S\to X$ is the inclusion and
$t:{\mathbb P}(N_{S|X}\oplus {\mathbb C})\to S$ is the projection.

Assume  $A _{+}=(A\cdot [E])F+\beta$ is an effective curve class of
${\mathbb P}(N_{S|X}\oplus {\mathbb C})$.
 Since the projection $t:{\mathbb P}(N_{S|X}\oplus {\mathbb C})\to S$ is
pseudo-holomorphic, $t_*((A\cdot [E])F+\beta)=t_*(\beta)$ is either
the zero class or an effective curve class of $S$ with respect to
$J|_S$. In the latter case, since $J|_S$ is compatible with
$\omega|_S$, $t_*\beta$ has positive symplectic area in $S$ with
respect to $\omega|_S$ and hence has positive symplectic area in $X$
with respect to $\omega$, i.e. $A _{-}=A-\beta$ is smaller than $A$.
Such a term on the right of (\ref{versus}) involves only standard
relative weighted graphs of $(\tilde X, E)$ lower in the partial
order than $ \Gamma( [pt], \varpi)|\mu$.


Let us focus on the terms with $t_*\beta=0$. In this case $\beta$
is a multiple of the fiber class $F$. Therefore $A _{+}$ is also
such a class.

Fix a splitting of $p_*A$ with $A _{+}$ a multiple of the fiber
class. Since the map $t\circ f _{+}: C _{+}\rightarrow S$ is
holomorphic, it maps every component of $C _{+}$ to either a point
or a holomorphic curve in $S$. Since they together represent the
zero class in $S$, each image must be a point. Hence the restriction
of $f _{+}$ to each connected component of $C _{+}$ also represents
a multiple of the fiber class.

  Next we show that $C _{-}$ must be a
connected curve of genus $g$. Since $C _{-}$ is a disjoint union of
connected components which forms a part of a degenerate genus $g$
Riemann surface, the sum of the arithmetic genera will be less than
or equal to $g$. If $C _{-}$ has more than one connected component,
its arithmetic genus will hence be strictly smaller than $g$  and
therefore the graph is of lower order.

\noindent{\bf Step IV--$\varpi$ insertions}. In this short step we
deal with the distribution of $\varpi$ insertions. If any of the
insertions of $\varpi_2$ is not empty, then the weighted relative
graph  of $(\tilde X,E)$ on the right of (\ref{versus}) is smaller
than  $ \Gamma( [pt], \varpi)|\mu$.

\noindent{\bf Step V--$\mu$ insertions}. In this step we show that
 $\eta$ in the largest weighted relative  graph with possibly nonzero coefficient is
 equal to $\mu$.
This step is a bit long and so we break into 4 sub-steps.

\noindent{\bf V.1 Set up}.
 It remains to analyze the
case in which the only non-relative insertions on the $\mathbb
P^k-$bundle side are given by $\tilde \mu$. From now on we focus
on the projective bundle side.

The distribution of the $l(\mu)$ insertions of $\tilde \mu$ among
the $l(\eta)$ rational components of $C _{+}$ decomposes the
relative insertion $\mu$ into $l(\eta)$ cohomology weighted
partitions
$$\pi^{(1)}, \cdots ,
\pi^{(l(\eta))},$$ with
$$\pi^{(r)}=\{(\mu_{n_1}, \delta_{i_1^{(r)}}),  \cdots , (\mu_{n_s},
\delta_{i_s^{(r)}})\},$$ and is allowed to be the empty partition.
Let $\tilde \pi^{(r)}$ be the corresponding absolute insertions.

Remember that there is only one relative marked point on each
component.
 Let $(\eta_r, \rho_r)$ be the
parts of $\eta$.

\noindent{\bf V.2 $S-$degree}. In this substep, we will show that
the largest term in (\ref{versus}) must have $\deg_S(\pi^{(r)}) +
\deg_S(\breve \rho_r) = \dim S$ for each $r$.

Notice that the curves in the fiber class live in the ${\mathbb
P}^k$ fibers.
 Then, for each $r$, in order for the multiple fiber class
relative invariant of the $\mathbb P^k-$bundle with the insertions
$<\tilde \pi^{(r)}| (\eta_r, \breve \rho_r)>$ to be nonzero, the
projections of cycles representing the Poincar\'e dual of $\delta$'s
and $\breve \rho_r$ have to intersect in $S$. This is to say that
$$\deg_S(\pi^{(r)})=\sum \deg_S(\delta_{i_j^{(r)}})\leq
\dim (S)-\deg_S(\breve \rho_r)=\deg_S(\rho_r).$$ By summing over all
$r$, we conclude that  $\deg_S(\mu)\leq \deg_S(\eta)$, and equality
holds if and only if $\deg_S(\pi^{(r)})=\deg_S(\rho_r)$ for all $r$.

Thus if $\deg _S(\rho_r)>\deg_S(\pi^{(r)})$ for some $r$, then the
relative invariant of $(\tilde X,E)$ on the right of (\ref{versus})
is smaller than  $ \langle [pt], \varpi|\mu\rangle_{g,A}^{  \tilde
X, E}$. It remains to analyze the case in which
\begin{equation}\label{dimS}\deg_S(\rho_r)=\deg_S(\pi^{(r)})=\sum
\deg_S(\delta_{i_j^{(r)}}),
\end{equation}
 i.e.
$\deg_S(\pi^{(r)})+\deg_S(\breve \rho_r)=\dim(S)$ for all $r$.

\noindent{\bf V.3 Dimension formula}.
 Since the $r-$th component of $C _{+}$ is totally
ramified with order $\eta_r$, it represents the class  $\eta_rF$.
Hence
 the dimension of the
moduli space for the $\eta_r-$totally ramified relative invariant of
the $\mathbb P^k-$bundle is
\begin{eqnarray*}
& &2<c_1({\mathbb P}^k\mbox{-bdl}), \eta_rF>+ \dim_{\mathbb R} (X)-6+2-2\eta_r+2l(\pi^{(r)})\\
&=&(2k+2)\eta_r+\dim_{\mathbb R}(X)-4-2\eta_r+2l(\pi^{(r)})\\
&=&2k\eta_r+\dim_{\mathbb R}(X)-4+2l(\pi^{(r)}).
\end{eqnarray*}
On the other hand, notice that
$$\deg(\tilde
\delta_i)=\deg_S(\delta_i)+2k.$$ Thus  the insertions $$<\tilde
\pi^{(r)}| (\eta_r,\breve \rho_r)>$$ have total degree
$$\begin{array}{ll}
&\deg(\tilde \pi^{(r)})+\deg(\breve \rho_r)+ \sum
_{j=1}^{l(\pi^{(r)})}2d_j(\mu)\\
=&2kl(\pi^{(r)})+\deg_S(\pi^{(r)})+\deg_S(\breve
\rho_r)+\deg_f(\breve \rho_r)\\ &+ \sum
_{j=1}^{l(\pi^{(r)})}2(k\mu_{n_j}^{(r)}-k+\frac{1}{2}\deg_f(\delta_{K_j}))\\
=&2\sum _{j=1}^{l(\pi^{(r)})}k\mu_{n_j}^{(r)} +
\dim(S)+\deg_f(\pi^{(r)})+\deg_f(\breve \rho_r).\\
\end{array}
$$
 Comparing the two formulas and using $\dim_{\mathbb R} X  =\dim_{\mathbb R} S +2k$, we conclude that
\begin{equation}\label{eta}
\begin{array}{ll}
&2k\eta_r+2k-4+2l(\pi^{(r)})\\
=&\sum _{j=1}^{l(\pi^{(r)})}2k\mu_{n_j}^{(r)}
 +\deg_f(\pi^{(r)})+\deg_f(\breve \rho_r).
\end{array}
\end{equation}

In fact, since we have assumed that (\ref{dimS}), we can arrive at
the same formula by simply computing the relevant relative invariant
of $({\mathbb P}^k, {\mathbb P}^{k-1})$.

Now the arguments depend on the value of $k$.

\noindent{\bf V.4 The case of  $k=1$}. For a ${\mathbb P}^1-$bundle,
we have $\deg_f=0, \deg_S=\deg$ and  $k=1$. The equation (\ref{eta})
simplifies to
\begin{equation} \label{eta_k}
\eta_r-1=\sum _{j=1}^{l(\pi^{(r)})}(\pi_j^{(r)}-1)= \sum
_{j=1}^{l(\pi^{(r)})}(\mu_{n_j}^{(r)}-1)
\end{equation}
 for each
$r$. Notice that when $\pi^{(r)}$ is empty, $l(\pi^{(r)})=0$ and the
right hand of formula (\ref{eta_k}) is understood to be 0.

Now consider the weighted partition $\pi^{(r)}$ containing a maximal
element $(\mu_1, \delta_{i_1})$ of $\mu$ in the size ordering.
Notice that $\mu_{n_j}^{(r)}\geq 1$ by the definition of relative
invariants. Hence,  by formula (\ref{eta_k}), either $\eta_r>\mu_1$,
or
$$\eta_r=\mu_1\hbox{ and all the other pairs of $\pi^{(r)}$ are of the
form $(1, \delta)$}.$$
 In the second case, according to
(\ref{dimS}), either $\deg(\rho_r)>\deg(\delta_{i_1})$, or
$$\deg
(\rho_r)=\deg(\delta_{i_1}), \hbox{ $\breve \rho_r$ is dual to
$\delta_{i_1}$},$$ and all other pairs of $\pi^{(r)}$ are of the
form $(1, {\bf 1})$.
 In fact, since the basis $\{\delta_i\}$ is
self-dual, we must have $\rho_r=\delta_{i_1}$.

Therefore either $\eta$ is larger than $\mu$ in the lexicographic
ordering and corresponds to a weighted relative graph of $(\tilde
X,E)$ strictly lower than  $ \Gamma( [pt], \varpi)|\mu$ in the
$\stackrel{\circ}{<}$ ordering, or the maximal pairs of $\eta$ and
$\mu$ agree.

 We now repeat the above analysis for the next largest
element of $\mu$ and continue until all the elements of $\mu$ not
equal to the smallest pair $(1, \delta_1)$ are exhausted.

Now let us understand how the smallest terms $(1, \delta_1)$ in the
lexicographic orderings in $\mu$ are distributed. We observe that
formula (\ref{eta_k}) sums to
\begin{equation} \label{eta-12}
\sum_r\eta_r-l(\eta)=\sum_j\mu_j-l(\mu).
\end{equation}
 Since $\sum_j \mu_j=A\cdot [E]$, we have
 $\sum_r \eta_r=\sum_j \mu_j$.
Any pair of the form $(1, \delta_1)$ also corresponds to some
$(\eta_r, \rho_r)$.
Thus $C _{+}$ has exactly $l(\eta)=l(\mu)$ connected components
 and  $\eta=\mu$, i.e. we recover $ \Gamma( [pt], \varpi)|\mu$ on the right of (\ref{versus}).


\noindent{\bf V.5 The case of  $k\geq 2$}. In this case $\pi^{(r)}$
cannot be empty by (\ref{eta}). This is because $\eta_r\geq 1$ and
hence $2k\eta_r+2k-4\geq 4k-4$, while $\deg_f\leq 2k-2$.

 If $\pi^{(r)}$ contains only one pair $(\mu_p, \delta_{i_p})$, then we have
 $\deg_S(\rho_r)=\deg_S(\delta_{i_p})$ by (\ref{dimS}), and in fact
 $\breve \rho_r$ is $S-$dual to $(\delta_{i_p})$.
If we write $\rho_r=\pi^*\theta\cup [E]^{\deg_f (\rho_r)}$ and
$\delta_{i_p}=\pi^*\theta'\cup [E]^{\deg_f(\delta_{i_p})}$, what we
have shown is that \begin{equation}\label{Spart}
\theta=\theta'\end{equation} as $\theta$ and $\theta'$ are elements
of  the self-dual basis of $S$.

 In addition, by (\ref{eta}),
\begin{equation} \label{last}2k\eta_r+2k-2=2k\mu_p+\deg_f(\delta_{i_p})+
 \deg_f(\breve \rho_r).
 \end{equation}
 Since $$0\leq \deg_f(\delta_{i_p})+
 \deg_f(\breve \rho_r)\leq 2k-2+2k-2,$$
 we have
 $$-(2k-2)\leq 2k(\eta_r-\mu_p)\leq (2k-2).$$
 Therefore
 $$\eta_r=\mu_p.$$ Furthermore, by (\ref{last}) it implies
 $$\deg_f(\delta_{i_p})+
 \deg_f(\breve \rho_r)=2k-2.$$
 It follows that $\deg_f
(\rho_r)=\deg_f(\delta_{i_p})$. Therefore by (\ref{Spart})  we must
have
$$\rho_r=\delta_{i_p}.$$

Since $\pi^{(r)}$ cannot be empty, if each $\pi^{(r)}$ contains at
most one pair, we recover the weighted relative graph $ \Gamma(
[pt], \varpi)|\mu$ on the right of (\ref{versus}).

 We argue now that if  some $\pi^{(r)}$ contains more
than one pair, then the corresponding weighted relative graph of
$(\tilde X, E)$ has either larger maximal tangency or  has the same
maximal tangency but larger $\deg_f$.

Now consider the weighted partition $\pi^{(r)}$ containing a maximal
element $(\mu_1, \delta_{i_1})$ of $\mu$ in the size ordering.
 We rewrite (\ref{eta}) as
$$\begin{array}{ll}
&2k(\eta_r-\mu_1)\\
=&\sum _{j\geq 3}^{l(\pi^{(r)})}2(k\mu_{n_j}^{(r)}-1)
+(2k\mu_{n_2}-2k)
 +\deg_f(\pi^{(r)})+
 \deg_f(\breve \rho_r).
\end{array}
$$
 Each term is non-negative since $\mu_i\geq 1$ and $k\geq 1$. As
$k$ is assumed to be at least $2$, all the terms are zero only if
$$l(\pi^{(r)})=2,\quad  \mu_{n_2}=1, \quad \deg_f(\pi^{(r)})=\deg_f(\breve
\rho_r)=0.$$ In particular, $\deg_f(\delta_1)=0$. Therefore the
relative invariant with insertions $<pt, \varpi|(\eta_r, \rho_r)> $
on the right of (\ref{versus}) is smaller since
$$\deg_f(\rho_r)=k-1-\deg_f(\breve
\rho_r)=k-1>0=\deg_f(\delta_1).$$

\noindent{\bf Step VI--Conclusion}. In summary  we have shown so far
that
\begin{equation} \label{transformation}
\begin{array} {lll}
& &\langle \tilde{\Gamma}( [pt], \varpi, \tilde \mu)\rangle^{  X}\\
& = & C_0  \langle \Gamma([pt], \varpi)|\mu\rangle^{  \tilde X, E}\\
& &+\mbox{ lower order terms,}\\
\end{array}
\end{equation}
where $C_0\ne 0$.

Viewing the invariants as the coordinates of the respective $\mathbb
R^{I_{ pt}}$, then (\ref{transformation}) defines a transformation
$${\mathcal A_S}:\mathbb R^{I_{ pt}}_{\tilde X, E}\to \mathbb R^{I_{ pt}}_{X, S}.$$
 As
$I_{pt}$ is indexed by the partial order, this matrix is lower
triangular with the $C_0$ as the diagonal entries. As the partial
ordered is lower bounded, this lower triangular matrix is
invertible.

\end{proof}

\begin{remark} When $k=1$ and $\sum_j \mu_j< A\cdot [E]$, the relative insertions of the largest
invariant on the right of (\ref{transformation}) are $\mu$ followed
by
 $A\cdot [E] -\sum_j \mu_j$ pairs of $(1, {\bf 1})$.  Notice that when $\sum_j \mu_j< A\cdot [E]$,  the relative
invariant $ \langle [pt], \varpi|\mu\rangle_{g,A}^{ \tilde X, E}$ is
zero by definition. What Theorem \ref{generalization} says in this
case is that $ \langle [pt], \varpi, \tilde \mu\rangle_{g,p_*(A)}^{
X}$ is expressed as the sum of standard relative invariants whose
weighted graph is lower than $\langle [pt],
\varpi|\mu\rangle_{g,A}^{ \tilde X, E}$.
\end{remark}

\section{Birational invariance}

In this section let $\tilde X$ be the blow up of $X$ along a
symplectic submanifold $Y$.
 Since GW invariants are unchanged under deformation, our main
theorem is equivalent to the following theorem.
\begin{theorem}  $\tilde X$ is uniruled if and only if $X$ is uniruled.
\end{theorem}

\begin{proof} Suppose $X$ is uniruled. Then
 there is a
 nonzero homology class $A\in H_2(X; {\mathbb Z})$, together with cohomology
classes $\alpha_2,  \cdots , \alpha_p\in H^*(X;\mathbb R)$ such that
the connected invariant
$$
     \langle [pt], \alpha_2, \cdots ,\alpha_p\rangle^{X}_{0,A} \not= 0.
$$
By the multi-linearity of GW-invariants of $X$,  we can assume that
$\alpha_i$ is one of the basis elements  $\sigma_{j_i}$ as in
(\ref{basisX}).
 Consider the degeneration
of $X$ into $(\tilde X, E)$ and ${\mathbb P}(N_{S|X}\oplus {\mathbb
C}), {\mathbb P}(N_{S|X}))$ and apply the degeneration formula of
3.3 to this invariant. If we put the point insertion on the $(\tilde
X, E)$ side and set
$$\alpha_i^-=\sigma_{j_i}^-=\gamma_{j_i}$$
 to be one of the basis elements in (\ref{basistildeX}), we find that
$\tilde X$ is uniruled relative to $E$. In fact, there is a nonzero,
possibly disconnected  genus 0 relative invariant of the form
\begin{equation}\label{nonzero relative invariant}
 \langle \Gamma^{\bullet}([pt], \varpi)|\mu\rangle^{\tilde X, E},
\end{equation}
where any insertion in $\varpi$ is of the form $\gamma_j$ with
$1\leq j \leq m_X$, $\mu$ is a standard weighted partition, and the
connected component containing the $[pt]$ insertion has nonzero
curve class. The relative invariant gives rise to a nonzero
component of the vector $v_{\tilde X, E}$ in ${\mathbb R}^{I_{0,
pt}}_{\tilde X, E}$. Apply the $I_{0,pt}$ version of Theorem
\ref{generalization} with $S=E$, together with the fact that the
blowup of $\tilde X$ along $E$ is still $\tilde X$ with a
deformation equivalent symplectic structure, we find that there is a
nonzero, possibly disconnected genus 0 absolute descendent invariant
of $\tilde X$ with a $[pt]$ insertion  in a connected component with
nonzero curve class. Hence $\tilde X$ is uniruled by Theorem
\ref{descendent}.

 Conversely, suppose that  $\tilde{X}$ is uniruled. Then there exists
a homology class $B\in H_2(\tilde{X};{\mathbb Z})$ and cohomology
classes $\beta_2,  \cdots, \beta_p\in H^*(\tilde{X};\mathbb R)$ such
that the connected invariant
$$
     \langle [pt], \beta_2, \cdots ,\beta_p\rangle^{\tilde{X}}_{0, B} \not= 0.
$$
 By the multi-linearity of GW invariants of
$\tilde X$, we can assume that $\beta_i$ is of the form
$\gamma_{j_i}$ as in (\ref{basistildeX}). Furthermore, assume that
$1\leq j_i\leq m_X$ if $2\leq i\leq l$, and $j_i\geq m_X+1$ if
$l+1\leq i \leq p$.

 Now apply  the degeneration formula to the degeneration of $\tilde
X$ along $E$ into $(\tilde X, E)$ and ($\mathbb P^1-$bundle, $E$),
distributing the first $[pt]$ insertion to the $(\tilde X, E)$ side
and the insertion $\beta_i$ with $i\geq l+1$ to the $\mathbb
P^1-$bundle side, i.e. setting $[pt]^+=0$ and  $\beta_i^-=0$ for
$i\geq l+1$. For  $\beta_i$ with $2\leq i\leq l$ we set
$\beta_i^-=\beta_i$. Notice that, for such a $\beta_i$, $\beta_i^+$
is not necessarily zero.
 Then
in the degeneration formula  there must be a nonzero relative
invariant of $(\tilde X, E)$ of the form (\ref{nonzero relative
invariant}) where  $\varpi$ has at most $l$ insertions of the form
$\gamma_j$ with $1\leq j \leq m_X$, $\mu$ is a standard weighted
partition, and the connected component containing the $[pt]$
insertion has nonzero curve class.

 Hence the
vector $v_{\tilde X, E}\in {\mathbb R}_{\tilde X, E}^{I_{0, pt}}$ is
nonzero. Now again apply the $I_{0,pt}$ version of Theorem
\ref{generalization} but this time with $S=Y$, we similarly conclude
that $X$ is uniruled.


\end{proof}

\begin{remark}
We could actually prove the following more precise result. If a
connected genus 0 relative invariant $ \langle \Gamma([pt],
\varpi)|\mu\rangle^{\tilde X, E}\ne 0$ and  has smallest  positive
symplectic area among all nonzero  invariants, then $ \langle
\tilde{\Gamma}([pt], \varpi,\tilde \mu)\rangle^{ \tilde X}$ is also
nonzero. Moreover, $A_{\Gamma}$ is the smallest uniruled class of
$\tilde X$.
\end{remark}

\begin{remark} It would be interesting to see whether the notion of
strongly uniruled in Remark \ref{strong} is a birational cobordism
property. There is also the notion of genus $g$ uniruled in
\cite{Lu2}. We could not prove this is a birational cobordism
property when $g$ is positive. What we can show is that the notion
of at most genus $g$ uniruled is a birational cobordism property.
\footnote{There is also the notion of strongly genus $g$ uniruled in
\cite{Lu2} and it is shown there that it is actually equivalent to
strongly uniruled.}
\end{remark}

\begin{remark}
 There is another important birational property of projective manifolds,
 the rational connectedness, which could be defined via the existence of
 a (connected, possibly reducible) rational curve  through any  two given points.
We could also define the symplectic analogue using connected genus
zero GW invariant with two point insertions. But it is not yet known
whether a projective rational connected manifold is symplectic
rational connected. Moreover, we could not show that this notion is
invariant under symplectic birational cobordisms. Our technique only
proves the disconnected version. But for rational connectedness,
unlike uniruledness, the disconnected version is strictly weaker.
\end{remark}

\section{Computing certain relative GW invariants of $({\mathbb P}^n, {\mathbb
P}^{n-1})$}

Suppose that $D$ is the hyperplane class of ${\mathbb P}^{n-1}$. In
this section, we compute the connected  genus zero relative GW
invariants of $({\mathbb P}^n, {\mathbb P}^{n-1})$,
\begin{equation} \label{pn}
\langle \tau_{nd-1-j}[pt]\mid D^j\rangle_{dL}^{{\mathbb
P}^n,{\mathbb P}^{n-1}},
\end{equation}
 which are factors of the diagonal entries
of ${\mathcal A}_S$.

\begin{theorem} \label{computation}
If $D$ is the hyperplane class of ${\mathbb P}^{n-1}$ and $L$ is
  the line class of ${\mathbb P}^{n}$
respectively, then for $d\geq 1$ and $0\leq j\leq n-1$,
$$
 \langle \tau_{nd-1-j}[pt]\mid D^j\rangle_{dL}^{{\mathbb
P}^n,{\mathbb P}^{n-1}} = \frac{1}{d^{n-j}[(d-1)!]^n}.
$$

\end{theorem}

 We now outline the proof of Theorem \ref{computation}.
 Choose homogeneous
coordinates
 $[z_0:z_1:\cdots:z_n]$ on ${\mathbb P}^n$. Let $(d)$ be the trivial partition
 of $d$ of length one. Denote by $H= \{ z_0=0\}$
the infinity hyperplane ${\mathbb P}^{n-1}_\infty$. Denote by
$\overline{\mathcal M}^{{\mathbb P}^n, H}_{0,1}(d, (d))$ the moduli
space of genus zero  relative stable maps to $({\mathbb P}^n, H)$
with one absolute marked point and one relative  marked point with
tangential order $(d)$. Here we simply choose the standard
integrable complex structure $J_{st}$ on $\mathbb P^n$. It is
well-known  that the moduli spaces of genus zero (absolute) stable
maps to $\mathbb P^n$ with this choice of almost complex structure
are actually smooth orbifolds of the correct dimension, i.e. they
represent the virtual moduli cycle (see e.g. \cite{FP}). Similar
arguments show that the genus zero moduli spaces of relative stable
maps to $(\mathbb P^n, H)$ with $J_{st}$ are smooth orbifolds of
corrected dimension as well. \footnote{We are indebted to R. Vakil
for confirmation on this.} In particular, $\overline{\mathcal
M}^{{\mathbb P}^n, H}_{0,1}(d, (d))$ is itself the virtual moduli
cycle  and hence
$$\langle\tau_{nd-1-j}[pt]\mid D^j\rangle_{dL}^{{\mathbb P}^n,
{\mathbb P}^{n-1}}\nonumber\\
  =  \int_{\overline{\mathcal M}^{{\mathbb P}^n, H}_{0,1}(d, (d))}  \psi^{nd-1-j}\wedge
 ev^{\mathbb P^n*}[pt]\wedge
ev^{H*}D^j\nonumber, $$
 where $ev^{\mathbb
P^n}$ is the evaluation map at the only absolute marked point $x_1$.

 Inspired by a calculation in the case $n=1$ in
\cite{OP}, we consider
$$V_1=\overline{\mathcal M}^{{\mathbb P}^n, H}_{0,1}(d, (d))\cap (  ev^{\mathbb P^n})^{-1}(p_0),$$
where $p_0=[1:0:\cdots :0]$.
 We then identify in Lemma \ref{x} a smooth orbifold $V_d$ in
$V_1$ and transform  the computation to an integral over $V_d$ in
Lemma \ref{appendix-1}. Finally we
 apply the (ordinary)
localization technique to evaluate  the integral in Lemma
\ref{appendix-2}.


Recall that $\mathcal L_1$ is the orbifold cotangent line bundle
whose fiber is the cotangent line at the unique absolute marked
point. We will simply  denote it by $\mathcal L$.

\begin{lemma} \label{x}
Let $V_i\subset V_1, 1\leq i\leq d $, be the subspace of relative
stable $J-$holomorphic maps with ramification order at least $i-1$
at the absolute marked point $x_1$. Then  $V_{i+1}\subset V_i$ is
the zero locus of a transverse section $s^{(i)}$ of $\oplus_n
{\mathcal L}^i$ over $V_i$, where $\oplus_n {\mathcal L}^{\otimes
i}$ is the direct sum of $n$ copies of the $i-$th power of $\mathcal
L$. In particular,  the smooth orbifold ${V_d}$ represents the class
\begin{equation}\label{Vd}
\frac{1}{[(d-1)!]^n}c_1({\mathcal L}_1)^{n(d-1)}
\end{equation}
in $V_1$.
\end{lemma}

\begin{proof}
 $\overline{\mathcal M}^{{\mathbb P}^n, H}_{0,1}(d,(d))$ consists of two
types of relative stable maps: the ones with the rigid target
${\mathbb P}^n$, and  the ones with a non-rigid target ${\mathbb
P}^n[m]$ for some $m\geq 1$. Recall that  ${\mathbb P}^n[m]$ is
obtained by gluing ${\mathbb P}^n$ with $m$ copies of the projective
bundle ${\mathbb P}({\mathcal O}_H(1)\oplus {\mathcal O})$ along the
zero sections and the infinity sections. Denote by $H_\infty$ the
last infinity divisor of ${\mathbb P}^n[m]$, which is just $H$ in
the case $m=0$.

  If $f\in V_d$ is mapped to a non-rigid target
  $\mathbb P^n[m]$ for some $m\geq 1$, due to the maximal ramification
  requirement at $x_1$ and $y_1$,
 the rubber part
must be a degree $d$ covering onto a chain of $\mathbb P^1-$ fibers
of ${\mathbb P}^n[m]$. But such maps  are invariant under
$\hbox{Aut}_H\mathbb P^n[m]\cong ({\mathbb C}^*)^m$ that dilates the
$\mathbb P^1-$fibers, and are thus not stable. So every relative
stable map in $V_d$ has rigid target ${\mathbb P}^n$, and in fact
$V_d$ is the $\mathbb Z_d-$orbifold parametrized by lines connecting
$p_0$ to points in $H=\mathbb P^{n-1}_{\infty}$.

For a relative stable map $ (C, x_1, y_1; f) \in V_1$ with $x_1$ and
$y_1$ as the absolute and relative marked points respectively,
 $x_1$ is mapped to $ p_0$ and the contact order of $f$ to $H_{\infty}$ at $y_1$ is $d$.

 Let $C_1$ be the irreducible component of $C$ containing $x_1$.
 Then the fiber $T^*_{x_1}C$ of the
orbifold complex line bundle ${\mathcal L}$ at $ (C, x_1, y_1; f)$
  is naturally identified with $T^*_{x_1}C_1$.
 If $C_1$ is contracted by
$f$ to $p_0$, then,

(i) $C_1$ must meet $\overline{C\setminus C_1}$ at no less than 2
points by stability, and

(ii) $\overline{C\setminus C_1}$ must be connected by the imposed
monodromy $(d)$ at $H_\infty$.

\noindent However  conditions (i) and (ii) violate the genus
constraint $g(C)= 0$. Thus the restriction of $f$ to $C_1$ is not a
constant map for any $f\in  V_1$.

 For $f\in V_1$ we have the pull-back map on the cotangent spaces:
\begin{equation}\label{appendix-0}
  f^*: T_{p_0}^*{\mathbb P}^n \longrightarrow T^*_{x_1}C_1=T^*_{x_1}C={\mathcal L}\mid_f.
\end{equation}
 Consider
homogeneous coordinates $[w_0:w_1]$ for $C_1\cong {\mathbb P}^1$
with $x_1 = [1:0]$. Take the complex coordinates
$[1:Z_1:\cdots:Z_n]$ where $Z_i = \frac{z_i}{z_0}$ over the
neighborhood $U_0 = \{z_0\not= 0\} \subset {\mathbb P}^n$ of $p_0$
and the complex coordinate $w= \frac{w_1}{w_0}$ over the
neighborhood $C_1- [0:1]$ of $x_1$. Then we have the canonical
isomorphisms,
$$
  T^*_{p_0}{\mathbb P}^n \cong \mathbb C^n= span_{\mathbb C}\{dZ_1, \cdots, dZ_n\}, \,\,\,\,\,\,\,\,
    T^*_{x_1}C_1  \cong span_{\mathbb C}\{dw\}.
$$
Write the restriction of $f$ to $C_1-[0:1]$ as $f(w) =
(Z_1(w),\cdots, Z_n(w))$. Observe that the map $f^*$ in
(\ref{appendix-0}) is the dual of the differential of $f(w)$ at
$w=0$. Thus the pullback map $f^*\in \hbox{Hom}(\mathbb C^n,
{\mathcal L}\mid_f)$ in (\ref{appendix-0}) yields the following
section $s^{(1)}$ of the direct sum $\oplus_n {\mathcal L}$ of $n$
copies of ${\mathcal L}$,
\begin{eqnarray*}\label{appendix}
   s^{(1)} : & V_1 & \longrightarrow \oplus_n{\mathcal L}\\
    & (C, x_1, y_1; f) & \longrightarrow (Z'_1(0)dw, \cdots,
   Z'_n(0)dw).
\end{eqnarray*}
Denote the zero locus $Z(s^{(1)})$ by $V_2$. Clearly $V_2\subset
V_1$ is the subspace of maps which have ramification order at least
$1$ at $x_1$.

Denote by $Ds^{(1)}$ the linearization of $s^{(1)}$, which  is
independent of the choice of local trivializations over
$V_2=Z(s^{(1)})$.

For $f\in V_1$ consider $n$ holomorphic tangent vector fields
$\xi_i\in T_fV_1, 1\leq i \leq n$ vanishing at $x_1$, and when
restricted to $C_1-[0:1]$, $\xi_i'(w) = (0,\cdots, 1,\cdots,0)$ with
the nonzero component only in the $i-$th entry.

Then we have $(Ds^{(1)})_f(\xi_i) = (0,\cdots, dw, \cdots, 0)$ with
the nonzero component only in the $i-$th   component. Thus
$Ds^{(1)}$ at a point of $Z(s^{(1)})$ is surjective. Therefore the
cycle $Z(s^{(1)})$ represents $c_1(\mathcal L)^{\otimes n}\cap
[V_1]$. When restricted to $Z(s^{(1)})$, the pull-back map on the
second differential of $f$ at $x_1$ will give rise to a canonical
section $s^{(2)}\in H^0(Z(s), \oplus_n {\mathcal L}^{\otimes 2})$,
i.e.
\begin{eqnarray*}
     s^{(2)}  : & V_2 &\longrightarrow \oplus_n {\mathcal L}^{\otimes
     2}\\
      & (C, x_1, y_1; f)&\longrightarrow (Z''_1(0)dw\otimes
     dw,\cdots, Z''_n(0) dw\otimes dw),
\end{eqnarray*}
where $\oplus_n {\mathcal L}^{\otimes 2}$ stands for the direct sum
of $n$ copies of ${\mathcal L}^{\otimes 2}$. It is easy to see that
$Z(s^{(2)}) \subset Z(s^{(1)})$ is the subspace where $x_1$ has
ramification order at least 2 over $p_0$. Hence the cycle
$Z(s^{(2)})$ represents the  class $[2 c_1({\mathcal L})]^n$.

After iterating the above construction, we conclude that
$[(d-1)!]^n\psi^{n(d-1)}$ is represented by $V_d$.

\end{proof}

Theorem \ref{computation} is an immediate consequence of Lemmas
\ref{appendix-1} and \ref{appendix-2}.

\begin{lemma}\label{appendix-1}
We have
\begin{equation*}
 \langle \tau_{nd-1-j}[pt]\mid D^j\rangle_{dL}^{{\mathbb P}^n,
{\mathbb P}^{n-1}} =
\frac{1}{[(d-1)!]^n}\int_{V_d}\psi^{n-1-j}\wedge ev^{H*}D^j,
\end{equation*} where $\psi=c_1(\mathcal L)$ and
$ev^{H*}$ is the evaluation map at the only relative marked point
$y_1$.
\end{lemma}

\begin{proof}
 By Lemma \ref{x} we have
\begin{eqnarray}\label{D'}
& & \langle\tau_{nd-1-j}[pt]\mid D^j\rangle_{dL}^{{\mathbb P}^n,
{\mathbb P}^{n-1}}\nonumber\\
 & = & \int_{\overline{\mathcal M}^{{\mathbb P}^n, H}_{0,1}(d, (d))}  \psi^{nd-1-j}\wedge
 ev^{\mathbb P^n*}[pt]\wedge
ev^{H*}D^j\nonumber \\
& = & \int_{V_1}\psi^{nd-1-j}\wedge ev^{H*}D^j \nonumber \\
 &=&
\frac{1}{[(d-1)!]^n}\int_{V_d}\psi^{n-1-j}\wedge
ev^{H*}D^j.\nonumber
\end{eqnarray}
\end{proof}

To evaluate the last  integral over the  smooth orbifold $V_d$ we
use localization.
\begin{lemma}\label{appendix-2} We have
$$\int_{V_d}\psi^{n-1-j}\wedge ev^{H*}D^j=\frac{1}{d^{n-j}}.$$
\end{lemma}

\begin{proof}
Consider the action of the group $T\cong (\mathbb C^*)^{n+1}$ of
diagonal matrices on ${\mathbb P}^n$ with pairwise distinct weights
$ \lambda_0, \cdots, \lambda_n$. The $T-$action on ${\mathbb P}^n$
has $n+1$ fixed points $p_i=[0:\cdots:0:1:0\cdots:0]$, $i=0, 1,
\cdots, n$. Denote by $\ell_{ij} = \ell_{ji}$, $i\not= j$, the line
in ${\mathbb P}^n$ passing through $p_i$ and $p_j$. Since the
infinity hyperplane $H$
 is invariant under the $T-$action, there is a natural $T-$action on
$\overline{\mathcal M}^{{\mathbb P}^n, H}_{0,1}(d, (d))$ by
translating the image of a relative stable map. The smooth orbifold
$V_d\subset \overline{\mathcal M}^{{\mathbb P}^n, H}_{0,1}(d, (d))$
is invariant under this  $T-$action, thus we may use the
localization technique to compute the integral. If we represent  the
Poincar\'e dual of $D^j$ by the $T-$invariant projective subspace of
${\mathbb P}^{n-1}$ generated by $p_1, \cdots, p_{n-j}$ (in some
sense localizing $D^j$ first), then we may replace $V_d$  by the
projective space ${\mathbb P}^{n-1-j}\subset V_d$ parameterizing
lines in ${\mathbb P}^n$ connecting $p_0$ and a point in ${\mathbb
P}^{n-1-j}$ generated by $p_1, p_2, \cdots, p_{n-j}$.

With the preceding understood, each fixed point in $V_d$
contributing to the integral corresponds to the following graph:

\begin{figure}[h]
\unitlength=1mm
\begin{picture}(50,10)(0,0)
\thicklines \put(6,5){\line(1,0){21}}
\put(6,5){\circle*{1.5}}\put(27,5){\circle*{1}}
\put(5,8){\{0\}}\put(15,7){d}\put(25,8){\{ 1\}}\put(25,1){i}
\end{picture}
\caption{Graph $\Gamma^i$ with $i\in \{1, 2, \cdots, n-j\}$}
\end{figure}
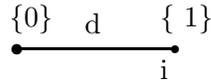
\noindent where the label $i$ indicates that the relative marked
point is mapped to the fixed point $p_i$. Therefore the normal
bundle to the fixed point corresponding to the graph $\Gamma^i$ is
the tangent space of ${\mathbb P}^{n-1-j}$ at $p_i$.  This  tangent
space  is generated by the tangent vectors of the line
$\ell_{i\alpha}$, $\alpha\in\{1, 2, \cdots, \hat{i}, \cdots, n-j\}$.
Then  the contribution of the virtual normal bundle to the fixed
point associated with the graph $\Gamma^i$  is
$$
  \prod_{\stackrel{\alpha=1}{\alpha\not=i}}^{n-j}\frac{1}{\lambda_i- \lambda_\alpha}.
$$
So by the Atiyah-Bott localization formula, we have
\begin{eqnarray}\label{appendix-3}
 \int_{V_d}\psi^{n-1-j}\wedge ev^{H*} D^j &=& \frac{1}{d}\sum_{i=1}^{n-j}\frac{(\frac{\lambda_i-\lambda_0}{d})^{n-1-j}}
{ \prod_{\stackrel{\alpha=1}{\alpha\not=i}}^{n-j}(\lambda_i- \lambda_\alpha)}\nonumber\\
& &\nonumber\\
&=&
\frac{1}{d^{n-j}}\sum_{i=1}^{n-j}\frac{(\lambda_i-\lambda_0)^{n-1-j}}
{ \prod_{\stackrel{\alpha=1}{\alpha\not=i}}^{n-j}(\lambda_i-
\lambda_\alpha)}=\frac{1}{d^{n-j}}.
\end{eqnarray}
The last equality follows from the expansion of the Vandermonde
determinant.
\end{proof}

 In a sequel paper we would need to compute some other relative invariants, and for
 that purpose we need to
apply the symplectic analogue of the virtual localization in
\cite{GV,GP}.

We end this section with the sketch of another argument of Lemma
\ref{appendix-2}. Let $V_{d, j}$ be the subspace of $V_d$ cut down
by $D^j$. Then $V_{d, j}$  can be identified with ${\mathbb
P}^{n-1-j}$, and  Lemma \ref{appendix-2} follows from the claim that
$\mathcal L^{\otimes d}$ is ${\mathcal O}(1)$ over $V_{d, j}$. To
prove this claim notice that $\mathcal L^{\otimes d}$ is the
cotangent line at $p_0$ of the ${\mathbb P}^1$ connecting $p_0$ and
a point $x$, which  corresponds to the plane ${\mathbb C}(p_0)\oplus
{\mathbb C}(x)$. And the tangent space at $p_0$ of this ${\mathbb
P}^1$ is the quotient of ${\mathbb C}(p_0)\oplus {\mathbb C}(x)$ by
${\mathbb C}(p_0)$, which is isomorphic to ${\mathbb C}(x)$. Thus
the dual of $\mathcal L^{\otimes d}$ is the tautological line bundle
${\mathcal O}(-1)$ on ${\mathbb P}^{n-1-j}$.

\end{document}